\documentclass[11pt,final]{siamltex}
\usepackage{graphicx}
\usepackage{float}
\usepackage{subfig}
\usepackage{amsmath}
\usepackage{amssymb}
\title{Alternative tilings for the fast multipole method on the plane}

\author{Yuancheng Luo and Ramani Duraiswami \thanks{Perceptual Interfaces and Reality Laboratory, Department of Computer Science, University of Maryland, College Park, MD, USA, \texttt{[yluo1,ramani]@umiacs.umd.edu}  } }

\begin{document}

\maketitle

\begin{abstract}
The fast multipole method (FMM) performs fast approximate kernel summation to a specified tolerance $\epsilon$ by using a hierarchical division of the domain, which groups source and receiver points into regions that satisfy local separation and the well-separated pair decomposition properties. While square tilings and quadtrees are commonly used in 2D, we investigate alternative tilings and associated spatial data structures: regular hexagons (septree) and triangles (triangle-quadtree). We show that both structures satisfy separation properties for the FMM and prove their theoretical error bounds and  computational costs. Empirical runtime and error analysis of our implementations are provided.
\end{abstract}



\pagestyle{myheadings}
\thispagestyle{plain}

\section{Introduction}
In the original FMM works by Greengard and Rokhlin \cite{FMM_GR}, multipole expansions are executed along hierarchical centers of hypercube arrangements that span an input domain. In the two dimensional case, the center of expansions can be mapped to a set of lattice points where their Vornonoi diagrams \cite{VORONOI} represent cell boundaries. That is, cell boundaries represent the equidistant points between nearest lattice points. For a square lattice configuration on the Euclidean plane, this equates to a square tiling and a general hypercube arrangement in higher dimensions. Composing or decomposing these tilings naturally leads to a self-replicating hierarchical quadtree data structure \cite{QUADTREE} that preserves local separation and well-separateness pair decomposition (WSPD) properties in \cite{WSPD}. This is necessary for the FMM to load-balance across levels and achieve linear runtime.

While the square lattice has been used by the FMM and the closely related treecode algorithm \cite{TREECODE} since their inceptions, we are not aware of other lattice configuration having been explored. In this paper, we investigate two alternative lattice groups that may be adapted into self-similar arrangements; Regular triangular and honeycomb lattice groups translate to regular hexagonal and triangle tilings. Although hexagonal tilings have been traditionally used in image processing \cite{HEXBOOK} and triangular tilings for mesh navigation \cite{TRIAG_QUADTREE}, they have not been applied to the FMM domain. We show how these arrangements form the bases for geometry preserving septree and triangle quadtree data structures that we adopt and implement for the hierarchical 2D FMM. Computational costs and error bounds induced by the data structures are derived to be comparable to that of quadtree.

The outline of this paper is as follows: Section \ref{sec:FMM} provides an algorithmic preface of a hierarchical FMM for two dimensional coulombic equations. Section \ref{sec:sepproperties} illustrates separation properties that emerge from the tilings. Section \ref{sec:datastruct} details a set of auxiliary functions shared amongst all data structures. Sections \ref{sec:septree} and \ref{sec:triagquad} introduce respective septree and triangle quadtree data structures and their implementations. Section \ref{sec:separationratios} derives level invariant separation ratios. Section \ref{sec:prop} provides an analysis of error bounds and computational costs. Section \ref{sec:experiments} presents empirical results that validate the theoretical analysis. Section \ref{sec:conclusions} concludes the paper and remarks on the generalizability of the data structures to higher dimensions.

\section{Fast Multipole Method} \label{sec:FMM}

The basic goal of the FMM is to evaluate the effects of a set of potentials denoted by source points $x_i$ with strengths $u_i$ on a set of receiver or target points $y_j$. A potential is approximated upto an error bound via a series of hierarchical multipole expansions arranged over the spatial domain. In section \ref{sec:experiments}, the algorithm is demonstrated on the 2D coulombic potential with kernel function
\begin{equation}
\begin{split}
\displaystyle
\Phi_{ji} = \log(y_j-x_i),
\end{split}
\end{equation}
over the complex plane though the discussion applies to general kernels. The total evaluation at a target point $j$ is 
\begin{equation}
\begin{split}
\displaystyle
\Phi_{j} = \sum_{i=1}^N \Phi_{ji} u_i, \quad j = \{ 1, \hdots, m \}, \quad y,x \in \mathbb{C}.
\end{split}
\end{equation}
An overview of the FMM algorithm is provided in \cite{FMM_COURSE}. For brevity, the expansion and translation operators in \cite{COULOMBIC_MAST} and \cite{COULOMBIC} for the coulombic kernel are omitted.

\section{Separation Properties} \label{sec:sepproperties}

\begin{figure}[ht]
\begin{center}
  \subfloat[Square]{\label{fig:square}\includegraphics[scale=0.2]{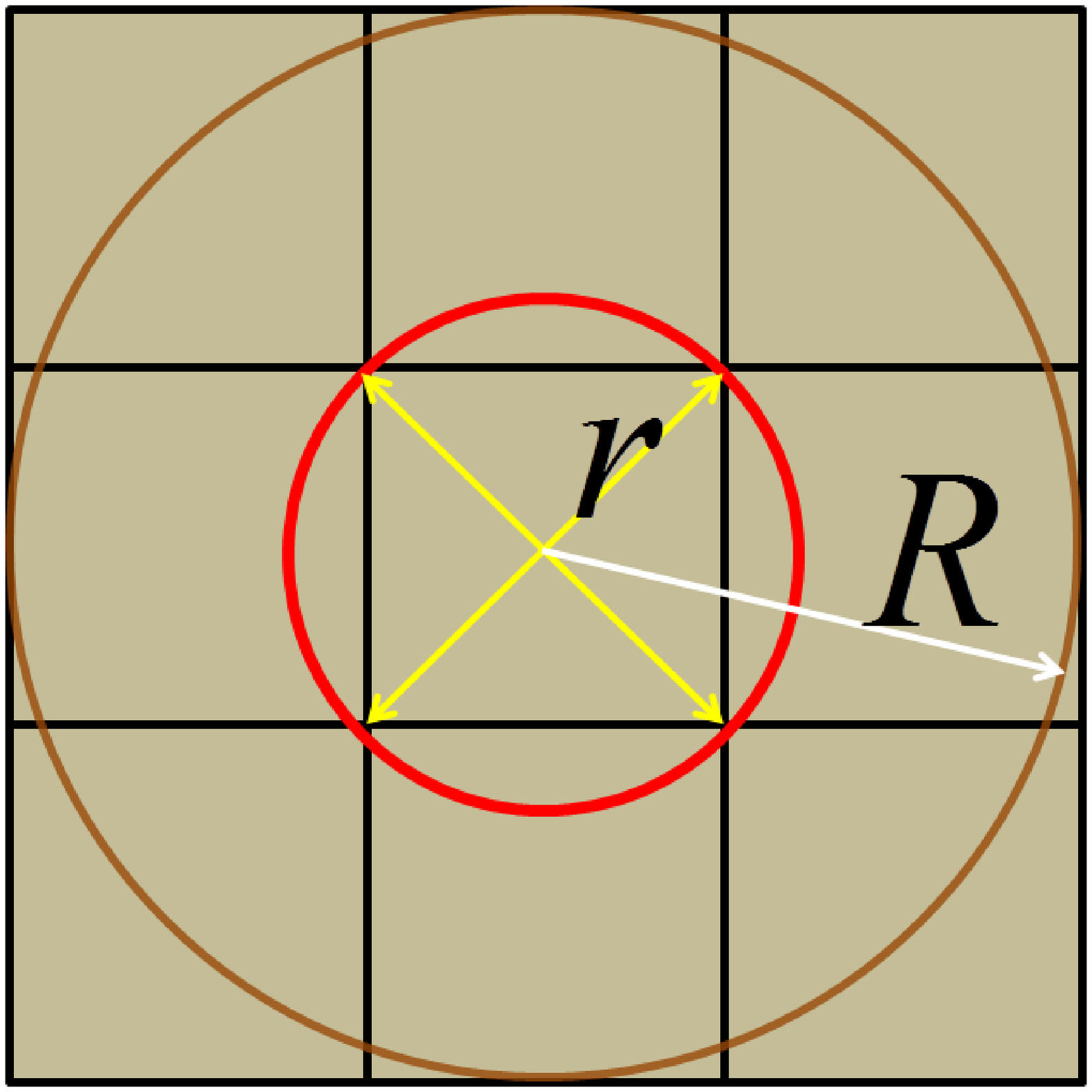}}    \hspace{.3cm}             
  \subfloat[Hexagon]{\label{fig:hexagon}\includegraphics[scale=0.2]{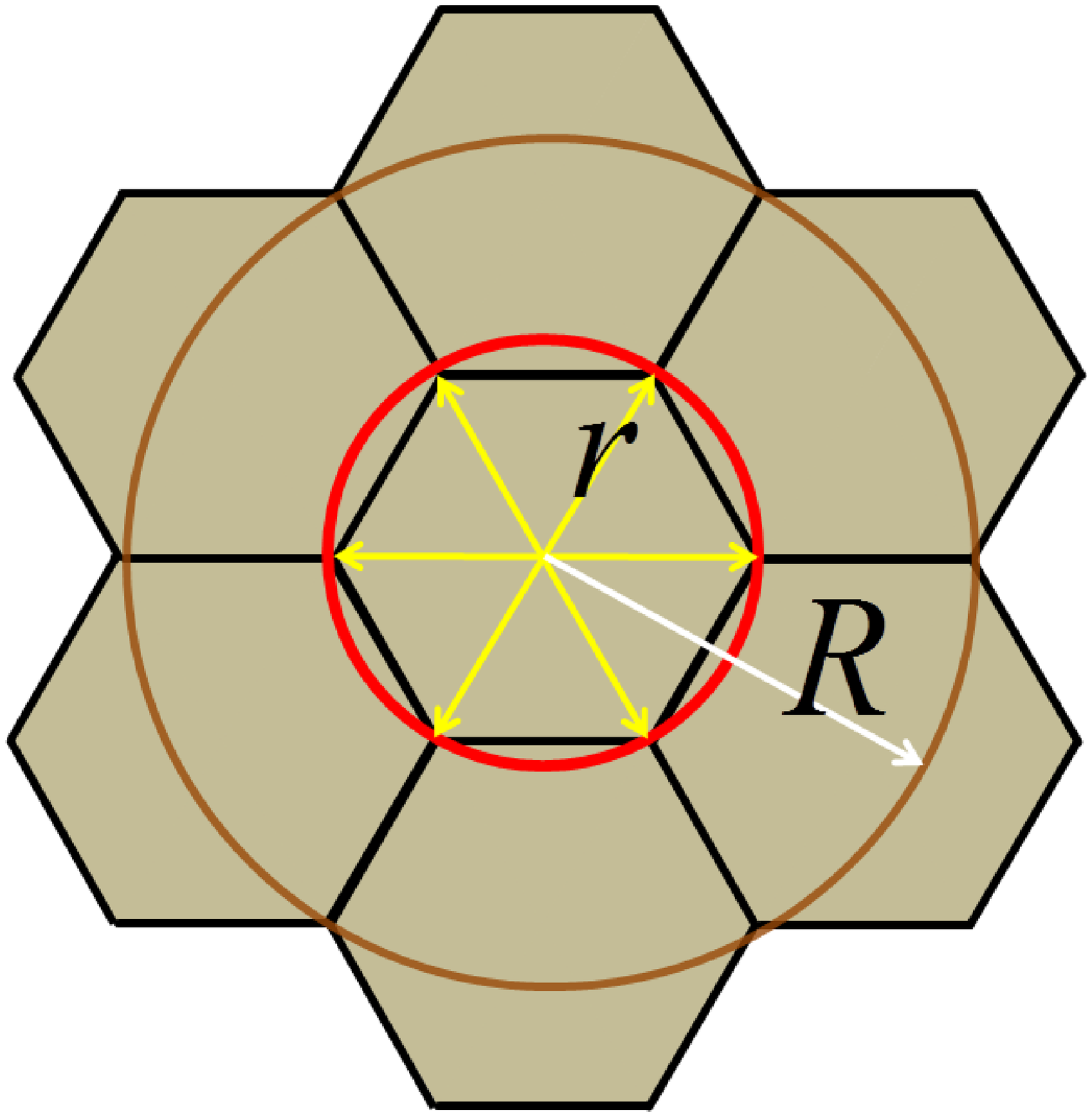}} \hspace{.3cm}
  \subfloat[Triangle]{\label{fig:triangle}\includegraphics[scale=0.2]{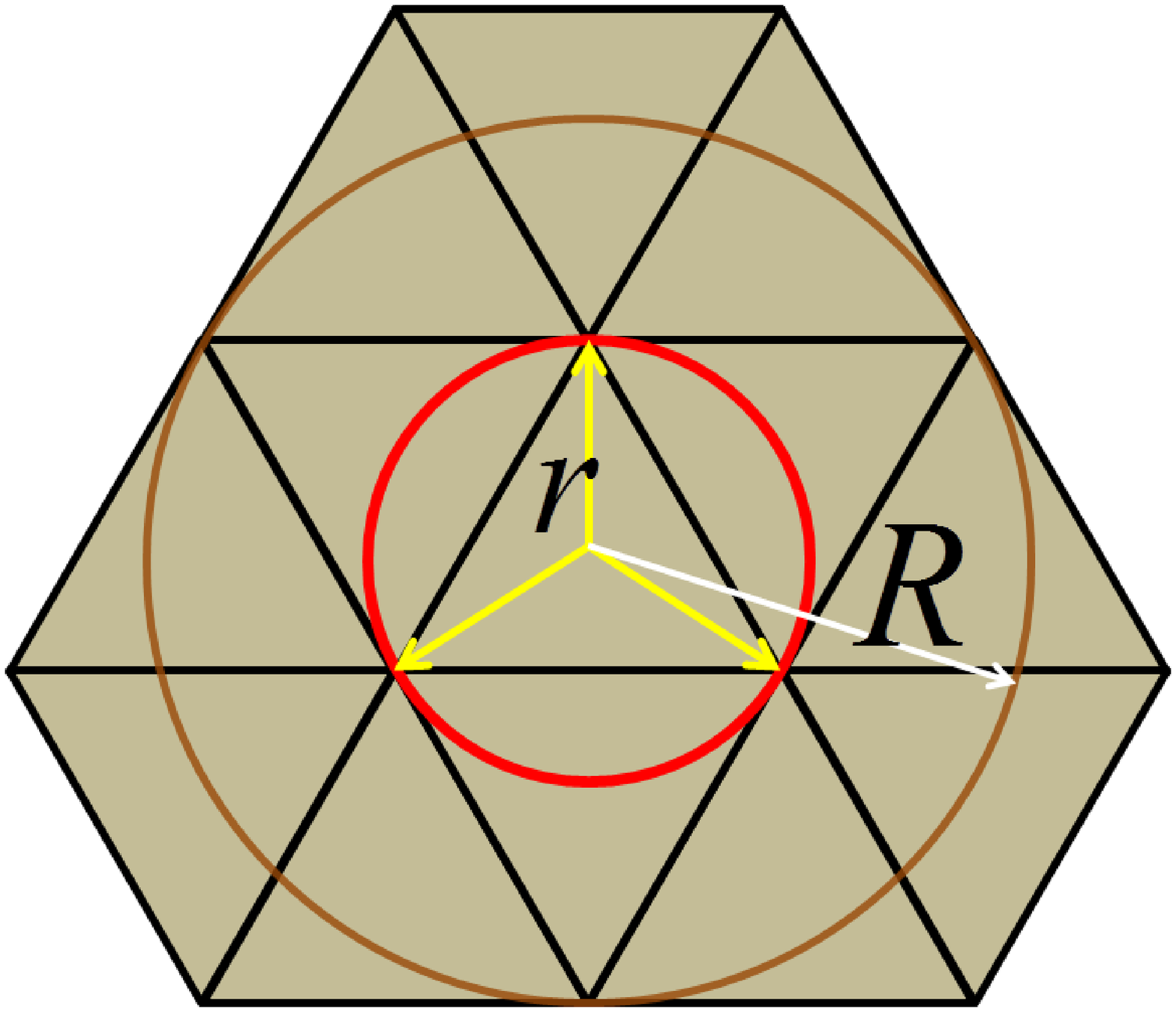}}
  \caption{Minor and major radii $r$ and major $R$ for the local separation property affects error in multipole expansions and translations}
  \label{fig:wellseparate}
  \end{center}
\end{figure}

A local separation property for multipole expansions and translations guarantees a minimum separation distance between points assigned to non-adjacent tiles. Geometrically, two circles that circumscribe and inscribe a tile and its adjacent neighbors in Figs. \ref{fig:wellseparate} bound a region. Formally, the separation ratio $r/R$ of radii between minor and major circles induces an error in the series expansion and so affects the number of truncation terms for a lower error bound provided in section \ref{sec:prop}. 

\begin{figure}[ht]
\begin{center}
\includegraphics[scale=0.25]{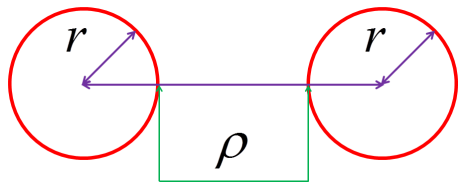}
\caption{WSPD distance $\rho$ between lattice points affects error in M2L translations}
\label{fig:rho}
\end{center}
\end{figure}

The WSPD property defines a minimum separation distance between multipole to local (M2L) centers for translation. Formally, the distance between separated sets in Fig. \ref{fig:rho} is expressed in terms of the ratio $\rho / r$ and used to estimate the number of truncation terms.

\section{Data Structures} \label{sec:datastruct}

To extend the local separation and WSPD properties beyond a basic tile, a concept of a cell is defined as a superset of tiles that shares a number of geometric and functional properties. For square and regular triangle tiles, cells are decomposable into smaller self-similar units. For square and regular hexagon tiles, cells are composable into larger aggregates. These hierarchical organizations give rise to indexable quadtree, septree, and triangle-quadtree data structures that when adapted for the FMM share the following attributes:
\begin{enumerate}
\item \emph{Separation}: Cells satisfy local separation and WSPD properties across all levels.
\item \emph{Indexing}: Cell indices satisfy some order relation in memory.
\item \emph{Spatial addressing}: Cell centers are computable from addresses and query points can find its bounding cell. 
\item \emph{Hierarchical addressing}: Cell children, parent, and vertex neighbor indices are computable from addresses.
\end{enumerate}

\begin{figure}[ht]
\begin{center}
\subfloat[\label{fig:quadtree}]{\includegraphics[scale=0.23]{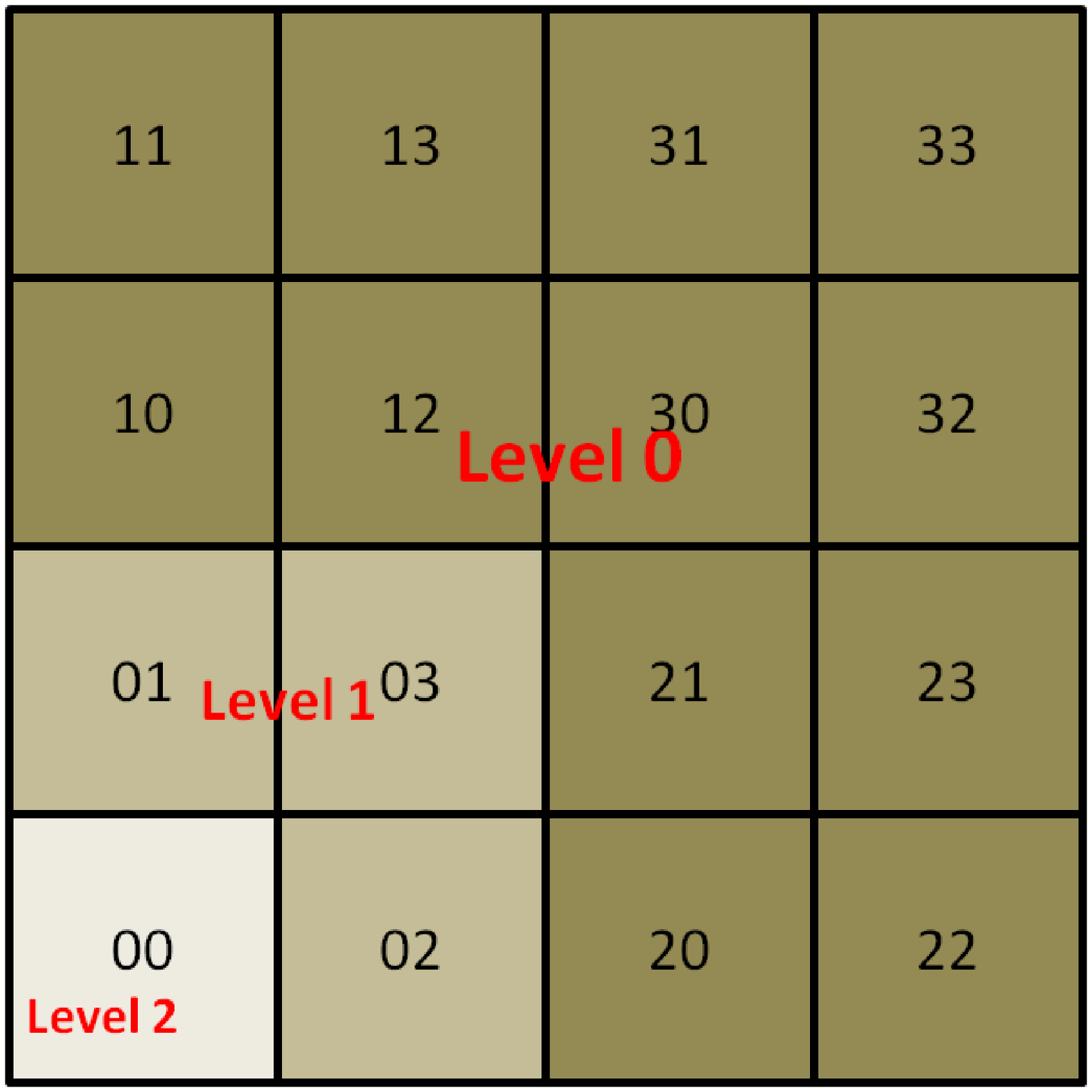}} \hspace{.3cm}
\subfloat[\label{fig:triagquad}]{\includegraphics[scale=0.23]{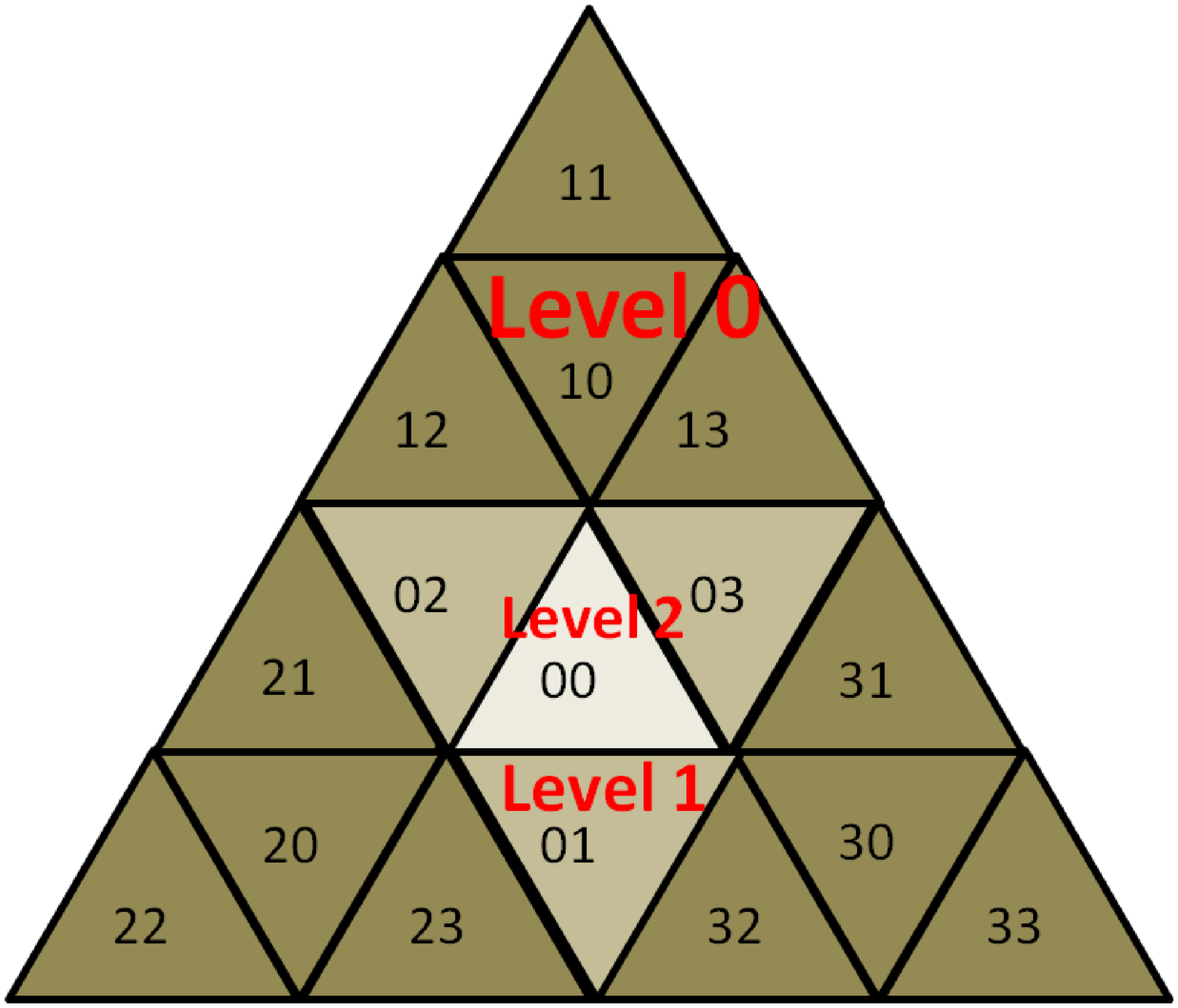}}
\caption{Morton $Z$ curve indexing for (a) quadtree  and (b) triangle quadtree}
\end{center}
\end{figure}

The separation properties for septree and triangle-quadtree are derived in section \ref{sec:separationratios}. The indexing property for the quadtree and triangle-quadtree follow Morton $Z$ curves in Figs. \ref{fig:quadtree} and \ref{fig:triagquad}. The septree follows a spiral pattern in Fig. \ref{fig:septree}.
The spatial and hierarchical addressing are accounted for by the following auxiliary functions:

\subsection*{CellCenter($n$, $l$)}
Returns the point coordinates of the center of cell $n$ on level $l$.

\subsection*{CellIndex($x$, $y$, $l$)}
Returns a cell $n$ on level $l$ that contains point $(x,y)$.

\subsection*{Children($n$)}
Returns cell indices of children cells of cell $n$.

\subsection*{Parent($n$)}
Returns a cell index of the parent of cell $n$.

\subsection*{Neighbors($n$, $l$)}
Returns cell indices on level $l$ that share a vertex with cell $n$ on level $l$.

\subsection*{NeighborsE4($n$, $l$)}
Finds cell indices required for M2L translation via
\begin{equation}
\begin{split}
\displaystyle
Children[Parent(n) \cup Neighbors(Parent(n), l-1)] \cap Neighbors(n, l).
\end{split}
\end{equation}

\section{Septree}\label{sec:septree}

\begin{figure}[ht]
\begin{center}
\subfloat[\label{fig:septree}]{\includegraphics[scale=0.33]{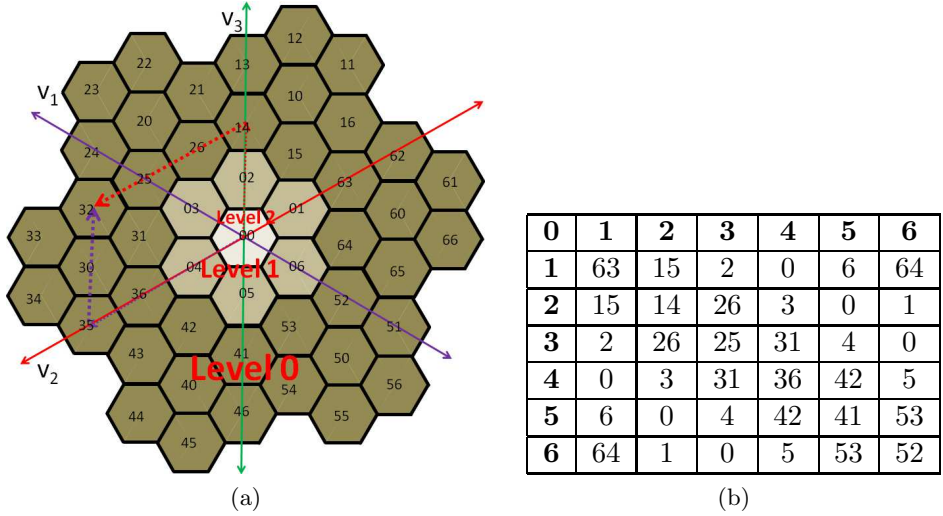}} \hspace{.3cm}
\subfloat[\label{tab:GBT}]{ \begin{tabular}[b]{ | c | c | c | c | c | c | c |} 
    \hline
    \textbf{0} & \textbf{1} & \textbf{2} & \textbf{3} & \textbf{4} & \textbf{5} & \textbf{6}\\ \hline
    \textbf{1} & 63 & 15 & 2 & 0 & 6 & 64\\ \hline
    \textbf{2} & 15 & 14 & 26 & 3 & 0 & 1\\ \hline
    \textbf{3} & 2 & 26 & 25 & 31 & 4 & 0 \\ \hline
    \textbf{4} & 0 & 3 & 31 & 36 & 42 & 5 \\ \hline
    \textbf{5} & 6 & 0 & 4 & 42 & 41 & 53 \\ \hline
    \textbf{6} & 64 & 1 & 0 & 5 & 53 & 52 \\ \hline
  \end{tabular}    }
  \caption{(a) Septree aggregates of hexagonal-like groups form a composable hierarchy where symmetries along the $v_1$, $v_2$, $v_3$ axis promote a base 7 GBT addressing scheme. (b) Table of all-pairwise GBT unit summations are analogous to vector summations to respective cells about the origin.}
  \end{center}
\end{figure}

The septree data structure is an upward composable hexagonal tessellation of the Euclidean plane. While the basic hexagon units cannot be subdivided into self-similar components, they may be aggregated into larger hexagonal-like groups. A base 7 indexing scheme begins on level $l_{max}$ with hexagon cell index $0_7$ centered at the origin and adjacent to neighbors with indices $\{ 1, \hdots, 6 \}$ tiled in a counter clock-wise direction. Aggregates of seven hexagonal cells on level $l_{max}$ form a single cell on level $l_{max}-1$ where they can be tiled and indexed in a similar fashion as seen in Fig. \ref{fig:septree}. Radial symmetry properties of the hexagon tessellation on level $l_{max}$ give rise to the base 7 Generalized Balanced Ternary (GBT) \cite{GBT} for cell indexing. This is investigated in section \ref{sec:septree:GBT} and its properties are exploited for a number of the FMM functions.

\subsection{Generalized Balanced Ternary} \label{sec:septree:GBT}
The generalization of Knuth's balanced ternary notation hierarchically describes permutohedral regions of $N-$dimensional spaces. In 2D, GBT indexing and arithmetic apply to hexagonal-like cells where the data can is decomposed. One important property realized by GBT addition is its duality with vector addition. That is, index arithmetic under GBT have a spatial analog with vector arithmetic in a coordinate system illustrated in Fig. \ref{fig:septree}.

For base 7 GBT unit addition, table \ref{tab:GBT} from Fig. \ref{fig:septree} represents all pair-wise summations of unit vectors taken in each of the six neighboring cell directions. Subsequent summations of higher order cell indices are interpreted as a mapping between coordinates along axes $v_1, v_2, v_3$ and their base 7 cell addresses. For notation, denote index $n_7$ as the base 7 expansion of the usual base 10 index $n$ and the GBT addition operator as $\circledast$, e.g. the sum of vectors to cell indices $14_7$ and $35_7$ along $v_2$ and $v_3$ axes is cell index $14_7 \circledast 35_7 = 32_7$.

\subsection*{Neighbors($n$, $l$)}
\begin{figure}[ht]
\begin{center}
\subfloat[\label{fig:septreeneigh}]{\includegraphics[scale=0.33]{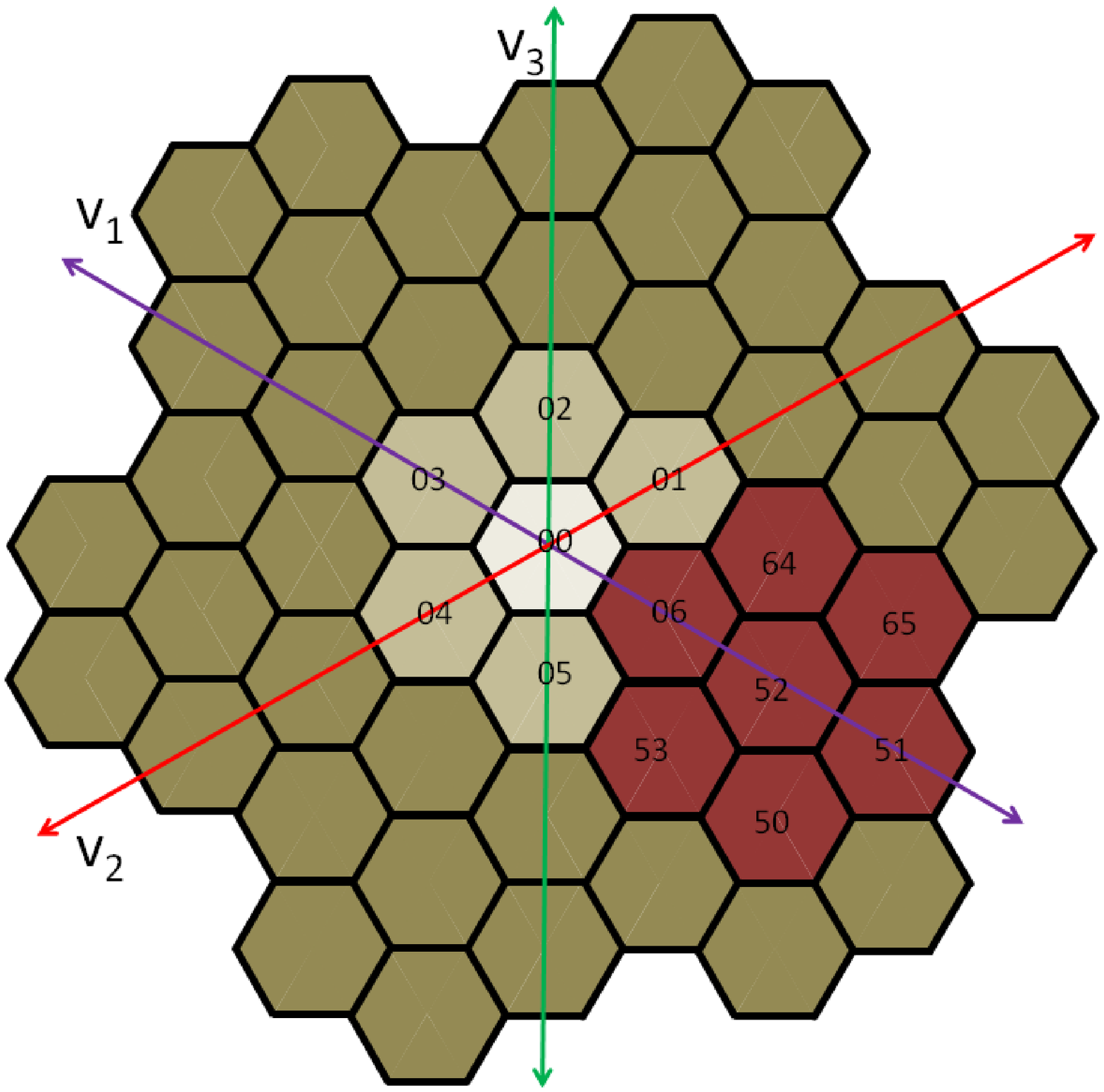}} \hspace{.3cm}
\subfloat[\label{fig:septreecenter}]{\includegraphics[scale=0.33]{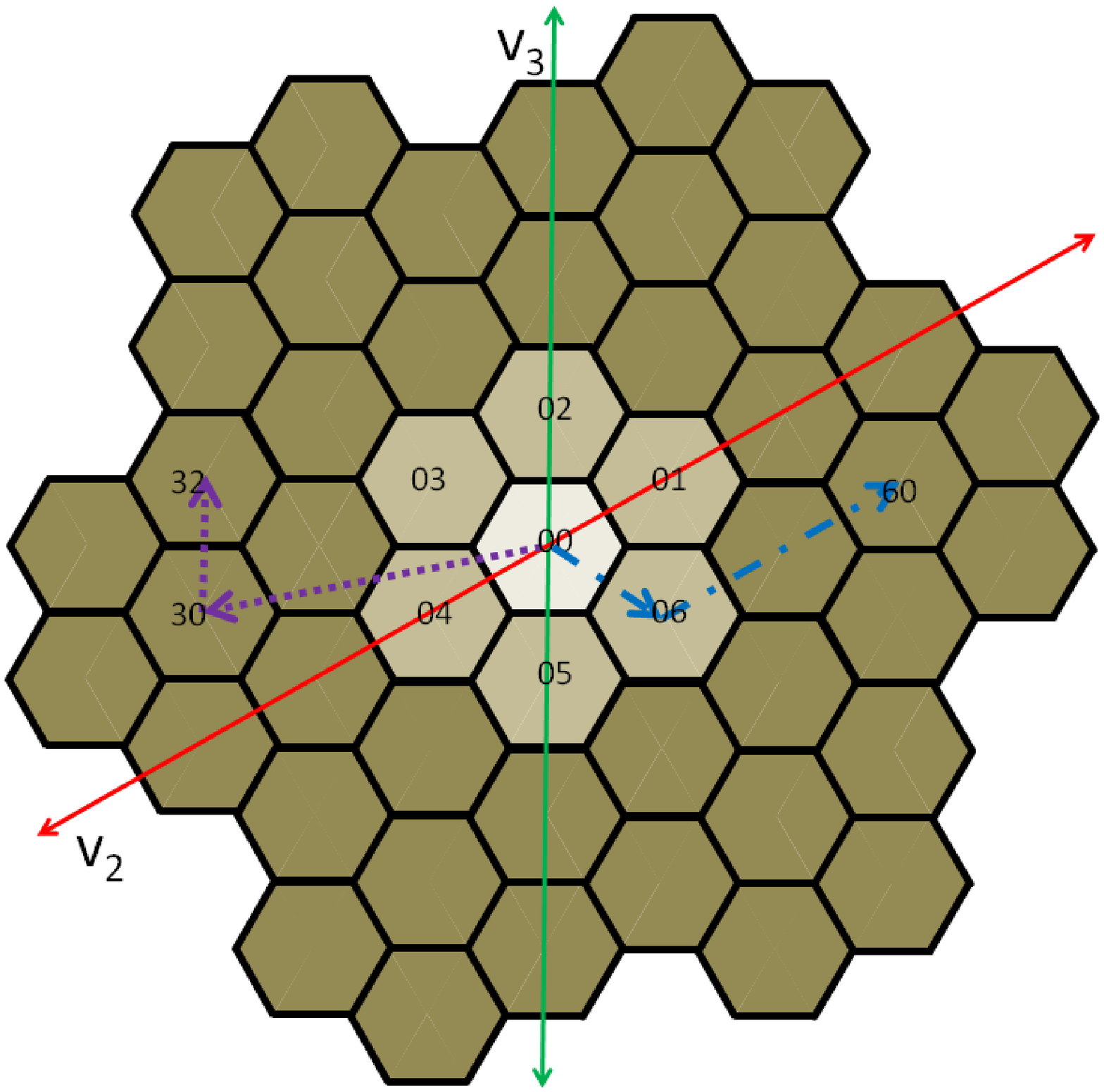}}
\caption{(a) Neighbors of cell index $52_7$ on level $2$ are $\{ 65, 64, 6, 53, 50, 51\}_7$. (b) Cell center of index $32_7$ found via translations along zero-extended components $30_2$ and $2_7$. Zero-extended index $60_7$ found via $6_7 \circledast  1_7 \circledast 1_7$}
\end{center}
\end{figure}

The level independent neighbors of cell $n$ are the pairwise GBT summations between cell index $n_7$ and the unit directions in Fig. \ref{fig:septreeneigh}. To adjust for the level, cell indices greater than $7^l-1$ are removed. The Neighbors($n$, $l$) procedure is as follows:
\begin{enumerate}
\item Let $N_7 = n_7 \circledast \{ 1, \hdots, 6\}$
\item Return indices in $N_{10}$ that are less than $7^l$
\end{enumerate}

\subsection*{CellCenter($n$, $l$)}
The center of cell $n$ on level $l$ is geometrically found via a series of translations along zero-extended components of index $n_7$. For notation, let index $t = n_7$ and component $t_1$ be the left-most digit. The cell center is expressed as
\begin{equation}
\displaystyle
\begin{split}
(x,y) = \sum_{i=1}^{|t|} T(t, i),
\end{split}
\end{equation}
where $|t|$ is total number of components and function $T:\mathbb{Z}^2 \rightarrow \mathbb{R}^2$ extends the $i^{th}$ component of $t$ with $|t|-i$ zeros and returns the corresponding vector to the cell from the origin, e.g. the cell center of index $342_7$ is the sum of vectors to the cell centers of indices $\{300, 40, 2 \}_7$.

\begin{figure}[ht]
\begin{center}
\subfloat[\label{tab:GBTextendedterms}]{ \begin{tabular}[b]{ |c| } 
    \hline
		$1 \circledast 2 \circledast 2 = 10$ \\ \hline
		$2 \circledast 3 \circledast 3 = 20$ \\ \hline
		$3 \circledast 4 \circledast 4 = 30$ \\ \hline
		$4 \circledast 5 \circledast 5 = 40$ \\ \hline
		$5 \circledast 6 \circledast 6 = 50$ \\ \hline
		$6 \circledast 1 \circledast 1 = 60$ \\ \hline
  \end{tabular}}\hspace{.3cm}
  \subfloat[\label{fig:GBTextendedgeometry}]{\includegraphics[scale=0.22]{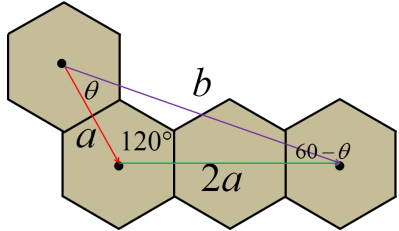}}
  \caption{(a) GBT unit addition for zero-extension. (b) Geometric analog for a zero-extension via translation $\vec{b}$}
  \end{center}
\end{figure}

The zero-extended $i^{th}$ component is the GBT summation of three shorter zero-extended components. The base case in table \ref{tab:GBTextendedterms} shows how successive GBT summations of a unit index $j$ and twice $1+(j \bmod{6})$ yield a zero-extended $j0$. For $i < |t|$, zero-extending the components ignores the preceding zeros, e.g. $10_7 \circledast 20_7 \circledast 20_7 = 100_7$. This is expressed by the relation
\begin{equation}
\begin{split}
\displaystyle
t_i 0^{|t|-i} = t_i 0^{|t|-i-1} \circledast [1+(t_i \bmod{6})]0^{|t|-i-1} \circledast [1+(t_i \bmod{6})]0^{|t|-i-1}.
\end{split}
\end{equation}

A geometric interpretation of index $t_i 0^{|t|-i}$ in Fig. \ref{fig:GBTextendedgeometry} reveals that the translation $\vec{b}$ is the sum of two vectors separated by $2 \pi /3$ radians with a two-to-one magnitude ratio. Solving for the angle of separation and the magnitude of vector $\vec{b}$ in appendix eqs. \ref{APP:GBTANG} and \ref{APP:GBTMAG} yield $\theta = \arctan{\frac{\sqrt{3}}{2}}$ and $b = a\sqrt{7}$. Hence, successive zero-extensions of a component $t_i$ are rotated by $\theta$ radians with magnitude $a_{i+1} =  a_i\sqrt{7}$ with a base magnitude $a_0 = r\sqrt{3}$.  

The function $T(t, i)$ can now be expressed in polar coordinates
\begin{equation}
\displaystyle
\begin{split}
j & = (|t|-i) + (l_{max} - l), \\
\theta_i & = j \arctan{\frac{\sqrt{3}}{2}} + t_i \frac{\pi}{3} - \frac{\pi}{6}, \quad R_i = r\sqrt{3} \left ( \sqrt{7} \right ) ^ j,
\end{split}
\label{eq:septreecenter}
\end{equation}
before $T(t,i) = (R_i \cos{\theta_i}, R_i \sin{\theta_i} ) $ gives the translation in a Cartesian coordinate system. The CellCenter($n$, $l$) procedure is as follows:
\begin{enumerate}
\item Let $(x,y) = (0,0)$,  $i=1$
\item Update $(x,y) = (x, y) + (R_i \cos{\theta_i}, R_i \sin{\theta_i} )$ from eq. \ref{eq:septreecenter} \label{algo:septreecenter:2} 
\item Increment $i$ and repeat from step \ref{algo:septreecenter:2} until $i=l$
\item Return point $(x,y)$
\end{enumerate}

\subsection*{CellIndex($x$, $y$, $l$)}
\begin{figure}[ht]
\begin{center}
\subfloat[\label{fig:septreeindex}]{\includegraphics[scale=0.33]{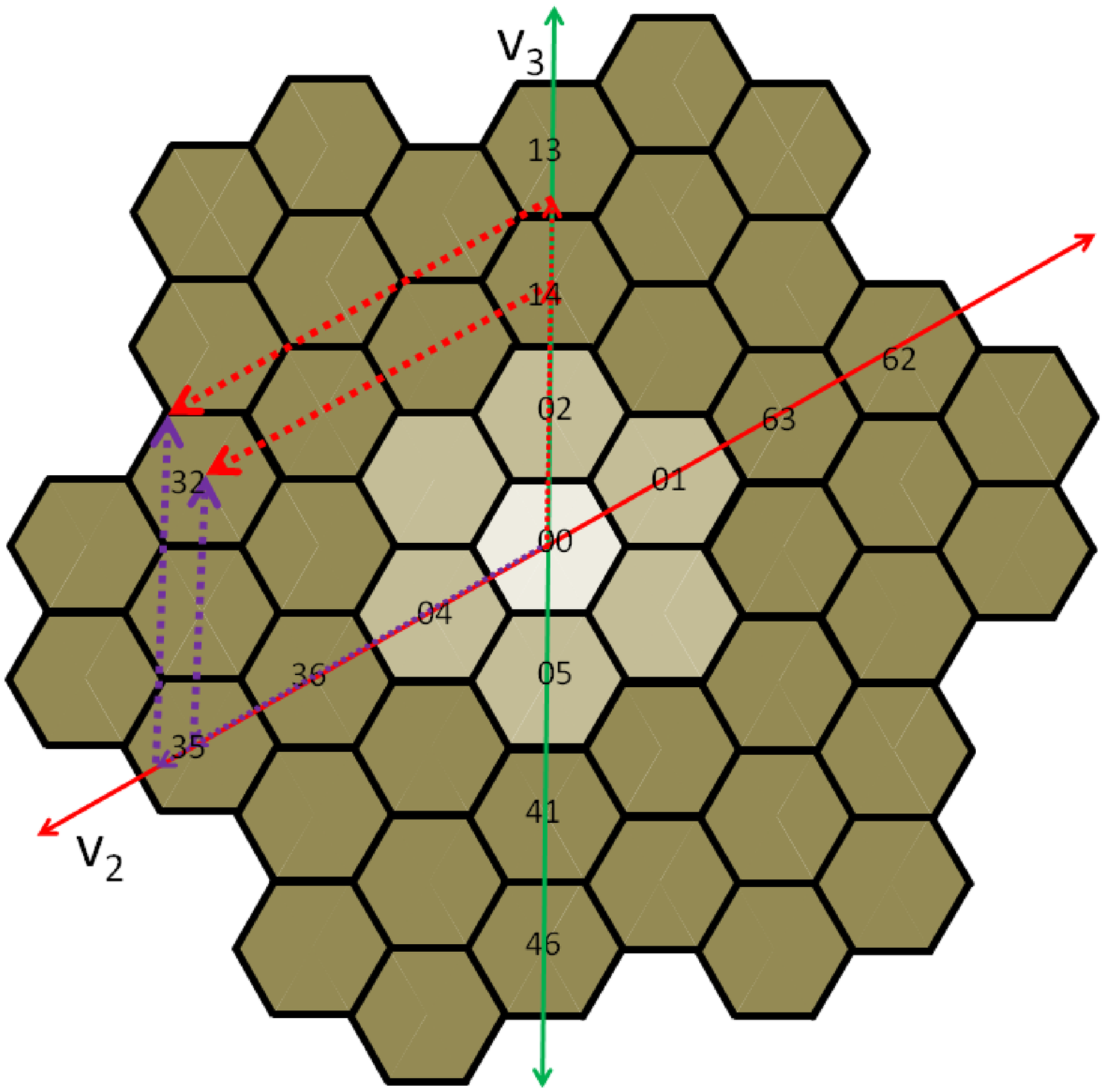}}\hspace{.3cm}
\subfloat[\label{fig:septreelattice}]{\includegraphics[scale=0.33]{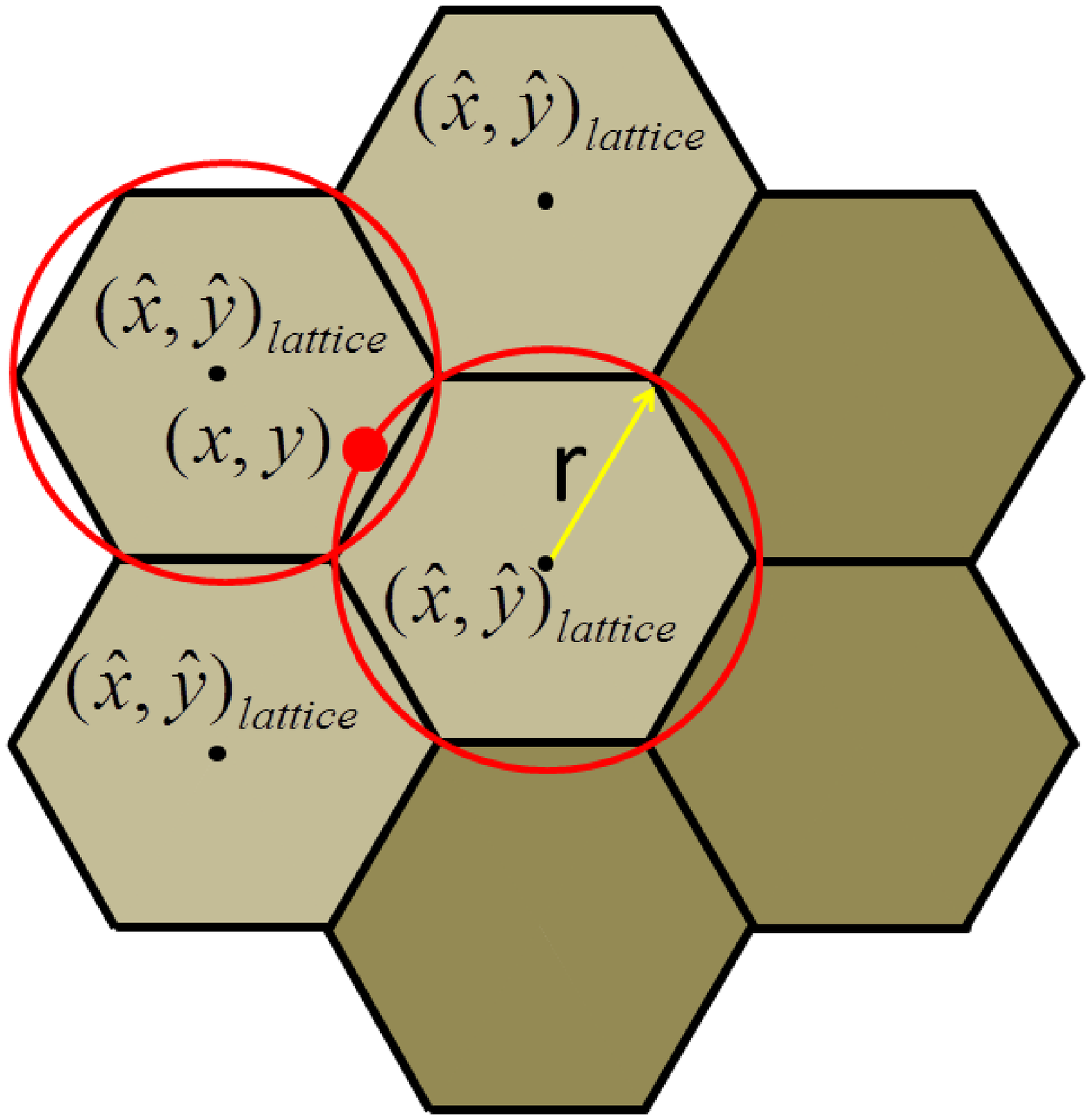}}
\caption{(a) Transforming query point $(x,y) \rightarrow (b,c)$ via change of basis may lead to incorrect cells along $v_2, v_3$ axes. (b) Hexagonal edges are Voronoi diagrams of equilateral triangular lattice points that construct cell centers}
\end{center}
\end{figure}

An approximation of the cell that contains the query point $(x,y)$ is found via a change of basis 
\begin{equation}
\displaystyle
\begin{split}
(b,c) = \left (  \frac{2x}{3r}, \frac{y\sqrt{3} - x}{3r}  \right ),
\end{split}
\label{eq:septreechangebasis}
\end{equation}
where coordinates $(b, c)$ are defined along axes $v_2, v_3$. These associated cells along the axes are erroneous as points near the corners of true bounding hexagon may project incorrectly in Fig. \ref{fig:septreeindex}.

To address this issue, boundary conditions are checked by observing that hexagonal edges are the Voronoi diagrams of equilateral triangular lattice points. That is, hexagonal edges represent points equidistant between nearest lattice points. The query point $(x,y)$ that falls within a hexagon would be the nearest neighbor of its center or lattice point. Furthermore, the query point may be bound between four neighbouring candidate lattice points in Fig. \ref{fig:septreelattice}. If the lattice points are defined w.r.t. axes $v_2, v_3$ as 
\begin{equation}
\displaystyle
\begin{split}
(\hat{x},\hat{y})_{lattice} = \left ( \frac{3r\hat{b}}{2} , \frac{r \sqrt{3}(2\hat{c}+\hat{b})}{2} \right ),
\end{split}
\end{equation}
then the candidate lattice points have coordinates 
\begin{equation}
\displaystyle
\begin{split}
\{ (\lfloor b \rfloor, \lfloor c \rfloor), (\lceil b \rceil, \lfloor c \rfloor), (\lfloor b \rfloor, \lceil c \rceil), (\lceil b \rceil, \lceil c \rceil) \}.
\end{split}
\label{eq:septreecanlattice}
\end{equation}

Finding the nearest lattice point w.r.t. the query point $(b,c)$ after converting back to Cartesian coordinates via eq. \ref{eq:septreechangebasis} yields the hexagon and its respective cell index $t = n_7$. The first $l$ components or $t_{1:l}$ is the the cell index on level $l$. The CellIndex($x$, $y$, $l$) procedure is as follows:
\begin{enumerate}
\item Compute axes coordinates $(b, c)$ from eq. \ref{eq:septreechangebasis}
\item Compute $i=1:4$ candidates lattice axes coordinates $(\hat{b}, \hat{c})_i$ from eq. \ref{eq:septreecanlattice}
\item Convert candidates lattice coordinates to Cartesian coordinates $(\hat{x}, \hat{y})_i$ with eq. \ref{eq:septreechangebasis}
\item Let the nearest candidate lattice point w.r.t. query point $(x,y)$ be $(\tilde{x}, \tilde{y})$ and equivalently $(\tilde{b}, \tilde{c})$
\item Find cell indices $u_7$, $v_7$ for coordinates $\tilde{b}$, $\tilde{c}$ along axes $v_2$, $v_3$ 
\item Return $u \circledast v$
\end{enumerate}

\section{Triangle Quadtree}\label{sec:triagquad}
A variation of the original quadtree data structure consists of a regular triangular tessellation of the Euclidean plane that is strictly decomposable. While each triangle may be subdivided into four similar units (center, vertical, left, right), a group of triangles is not upward composable in the sense that vertical orientations of descendant triangles depend on the orientation of their common ancestor triangle or root. That is, all children of a center cell will have inverted its orientation. For notation, denote the default upright orientation as $up_i=1$ and the inverted orientation as $up_i=-1$ for a triangle on level $i$. Triangle quadtree indexing begins on level $i=0$ with cell index $0$ centered at the origin and enclosing the entire domain. Children cells on level $i+1$ with indices $\{ 0, \hdots, 3\}$ are associated with cell types $\{\textrm{center, vertical, left, right} \}$ ordered w.r.t. the parent cell.

\subsection*{Neighbors($n$, $l$)}

\begin{figure}[ht]
\begin{center}
\subfloat[\label{fig:triagquadadjneigh}]{\includegraphics[scale=0.25]{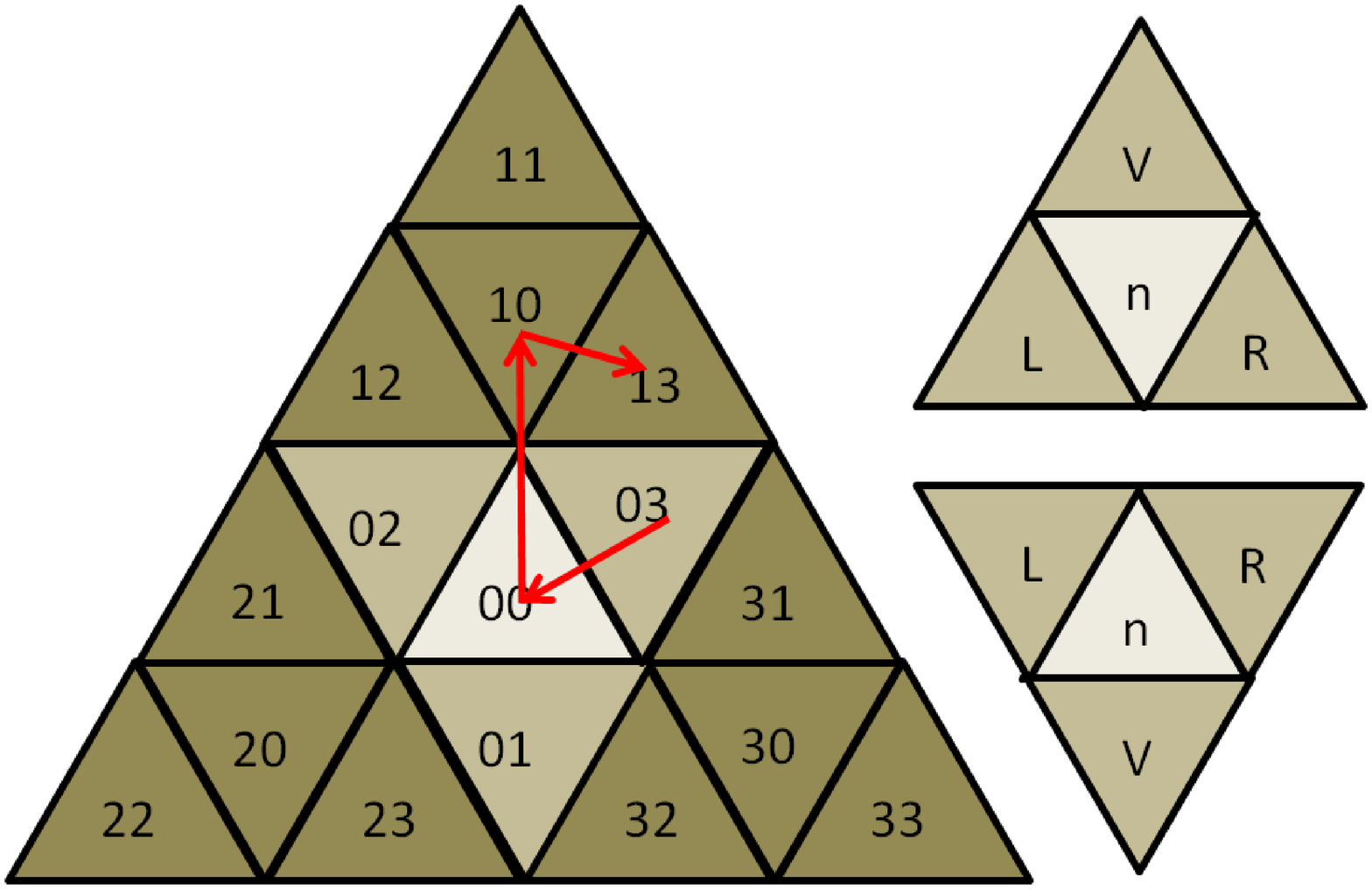}} \hspace{.3cm}
\subfloat[\label{fig:triagquadneigh}]{\includegraphics[scale=0.25]{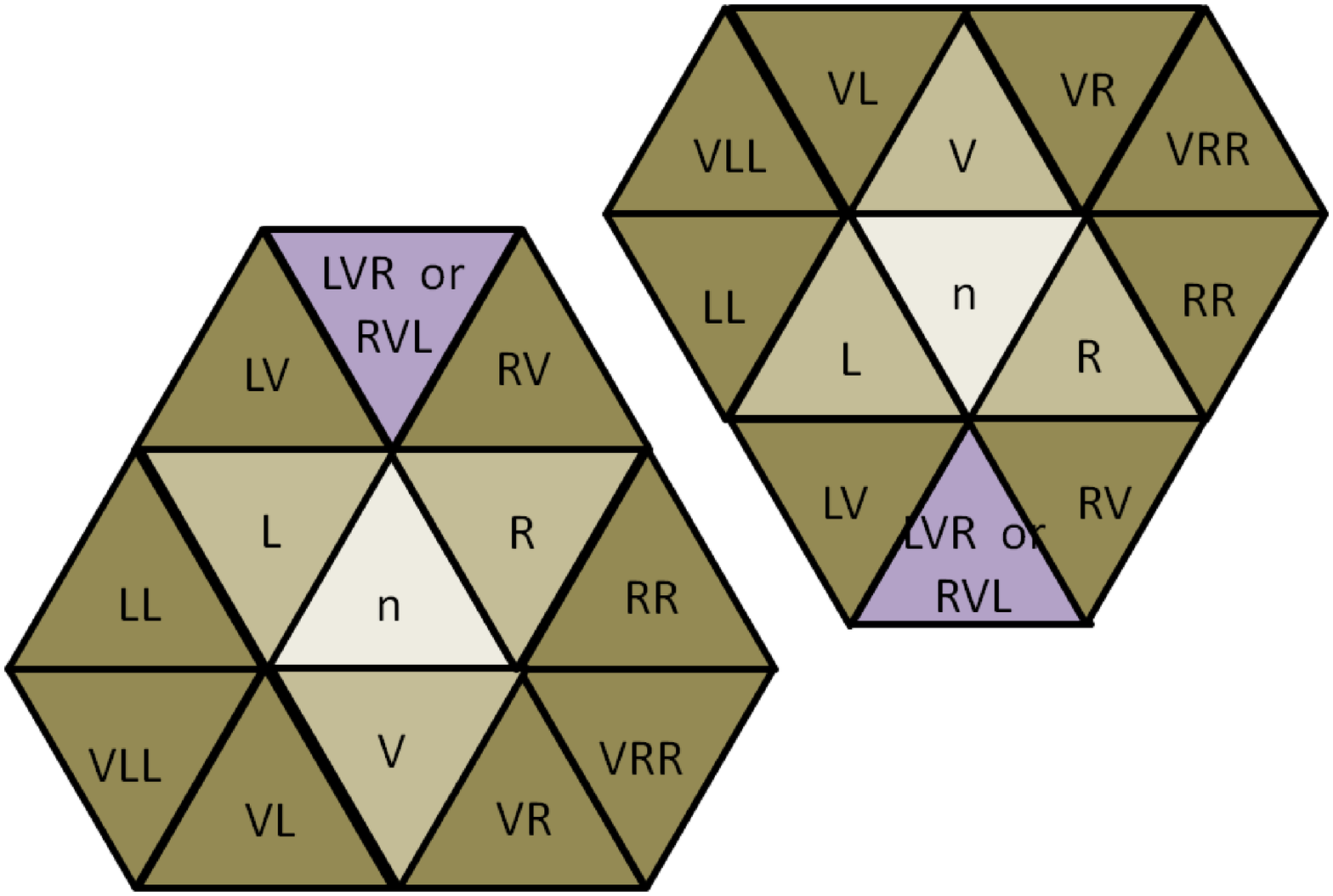}}
\caption{(a) Neighbor finding from cell index $03$ in the vertical direction. Cell $0$ is nearest ancestor with vertical sibling cell $1$ on level $1$. Reflect the path leading to ancestor cell $0$ gives cell index $13$. (b) Neighbor finding from recursive calls to adjacent neighbors}
\end{center}
\end{figure}

For regular quadtrees, neighbors of a cell index $n$ need only to share an edge. For triangle quadtrees, denote the left, right, and vertical adjacent neighbors as $ \{L, R, V \}$ in Fig. \ref{fig:triagquadadjneigh}.

Similar to \cite{TRIAG_QUADTREE}, the adjacent neighbor of cell $n$ in direction $k$ is found with an ancestor-sibling-reflect method in the adjacentNeighbor($n$, $k$) procedure:
\begin{enumerate}
\item Search for the closest common ancestor containing a sibling in direction $k$. For notation, let cell index $t = n_4$ and matrix $D$ in table \ref{tab:triagadjneisearch} map cell source types and search directions to destination cell types. From right to left components $t_i$, find the first valid entry $D_{t_i, k}$.

\begin{table}
\begin{center}
  \begin{tabular}{ | c | c | c | c |} 
    \hline
    \textbf{Source cell type} & \textbf{Left cell ID} & \textbf{Right cell ID}  & \textbf{Vert cell ID} \\ \hline
   0 & 2* & 3* & 1* \\ \hline
   1 & 3 & 2 & 0* \\ \hline
   2 & 1 & 0* & 2 \\ \hline
   3 & 0* & 1& 3 \\ \hline   
  \end{tabular}    
  \caption{Mapping from source cell types and search directions to destination cell types. * entries indicate valid siblings of a source cell}
    \label{tab:triagadjneisearch}
\end{center}
\end{table}

\item Move to sibling node. Replace component $t_i$ with entry $D_{t_i, k}$. 
\item Reflect the path taken to reach ancestor. For $j = \{ (i+1), \hdots, l_{max} \}$, replace component $t_j$ with entry $D_{t_j, k}$.
\item Return cell index $t$
\end{enumerate}

To find all neighbors that share a vertex with cell $n$, make a set of recursive calls to adjacentNeighbor($n$, $k$) in Fig. \ref{fig:triagquadneigh}. The Neighbors($n$, $l$) procedure is as follows:
\begin{enumerate}
\item Let $\{ L, R, V \}$ be adjacent neighbors of cell $n$ in directions left, right, vertical
\item Let $\{ VL, VR, LL, RR, LV, RV \}$ be adjacent neighbors of cells $V, L, R$ in directions left, right, vertical
\item Let $\{ VLL, VRR, LVR \}$ be adjacent neighbors of cells $VL, VR, LV$ in directions left and right
\item Return cell indices $\{ L, R, V, VL, VR, LL, RR, LV, RV, VLL, VRR, LVR \}$
\end{enumerate}

\subsection*{CellCenter($n$, $l$)}
The center of cell $n$ on level $l$ is found via a series of translations for each component of index $n_4$. Let index $t=n_4$, component $t_1$ be the left-most digit, and cell center 
\begin{equation}
\displaystyle
\begin{split}
(x,y) = \sum_{i=1}^{|t|} F(t, i),
\end{split}
\end{equation}
where function $F:\mathbb{Z}^2 \rightarrow \mathbb{R}^2$ returns a vector to cell index $t_i$ that is adjusted for level and orientation, e.g. the cell center of index $302_4$ is the sum of vectors to cell centers of indices $\{ 3, 0, 2 \}_4$ on levels $\{ 1, 2, 3\}$ with orientations $\{ 1, -1, -1\}$.

\begin{figure}[ht]
\begin{center}
\subfloat[\label{fig:triagquadcenter}]{\includegraphics[scale=0.33]{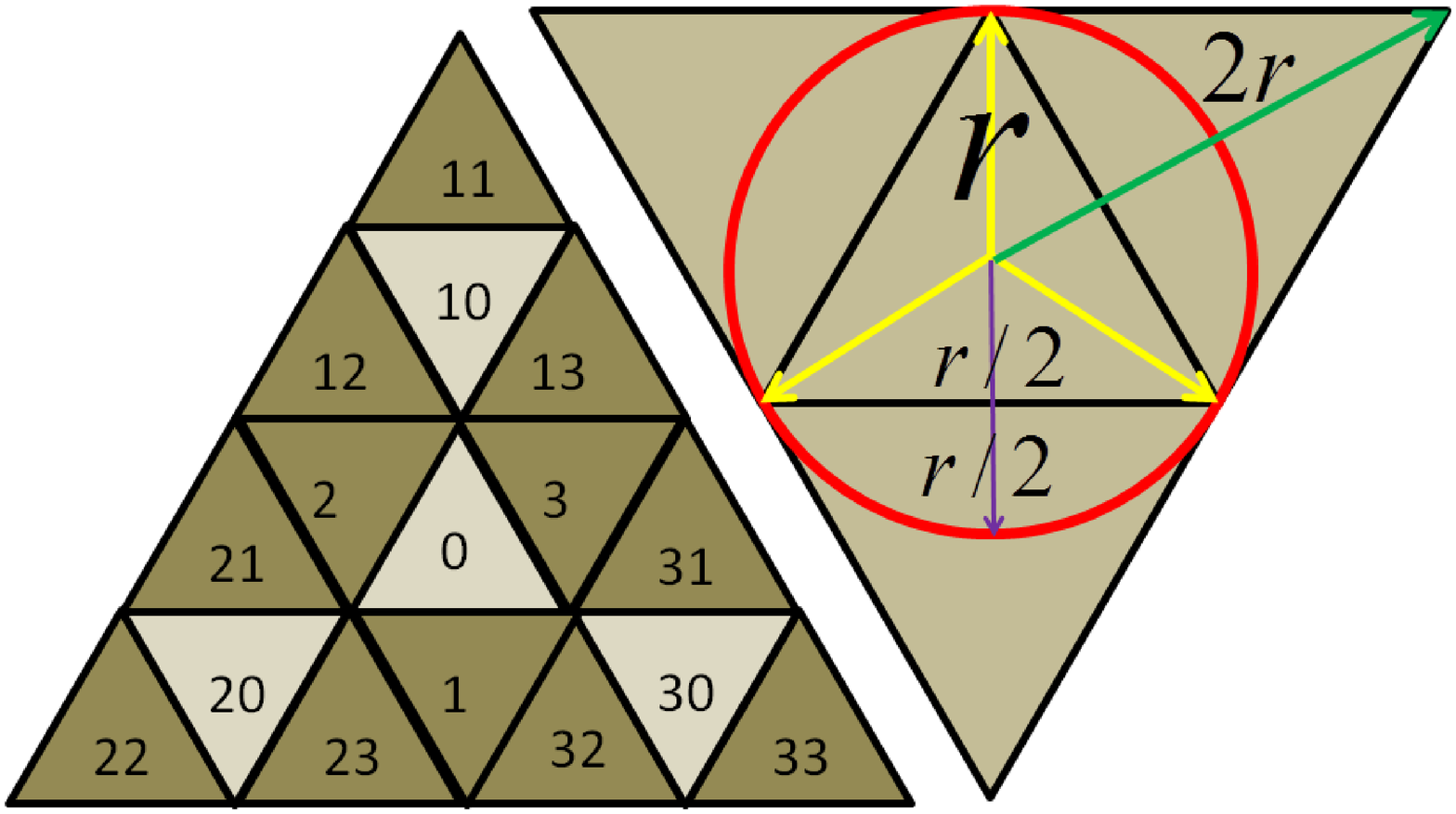}} \hspace{.3cm}
\subfloat[\label{fig:triagquadindex}]{\includegraphics[scale=0.3]{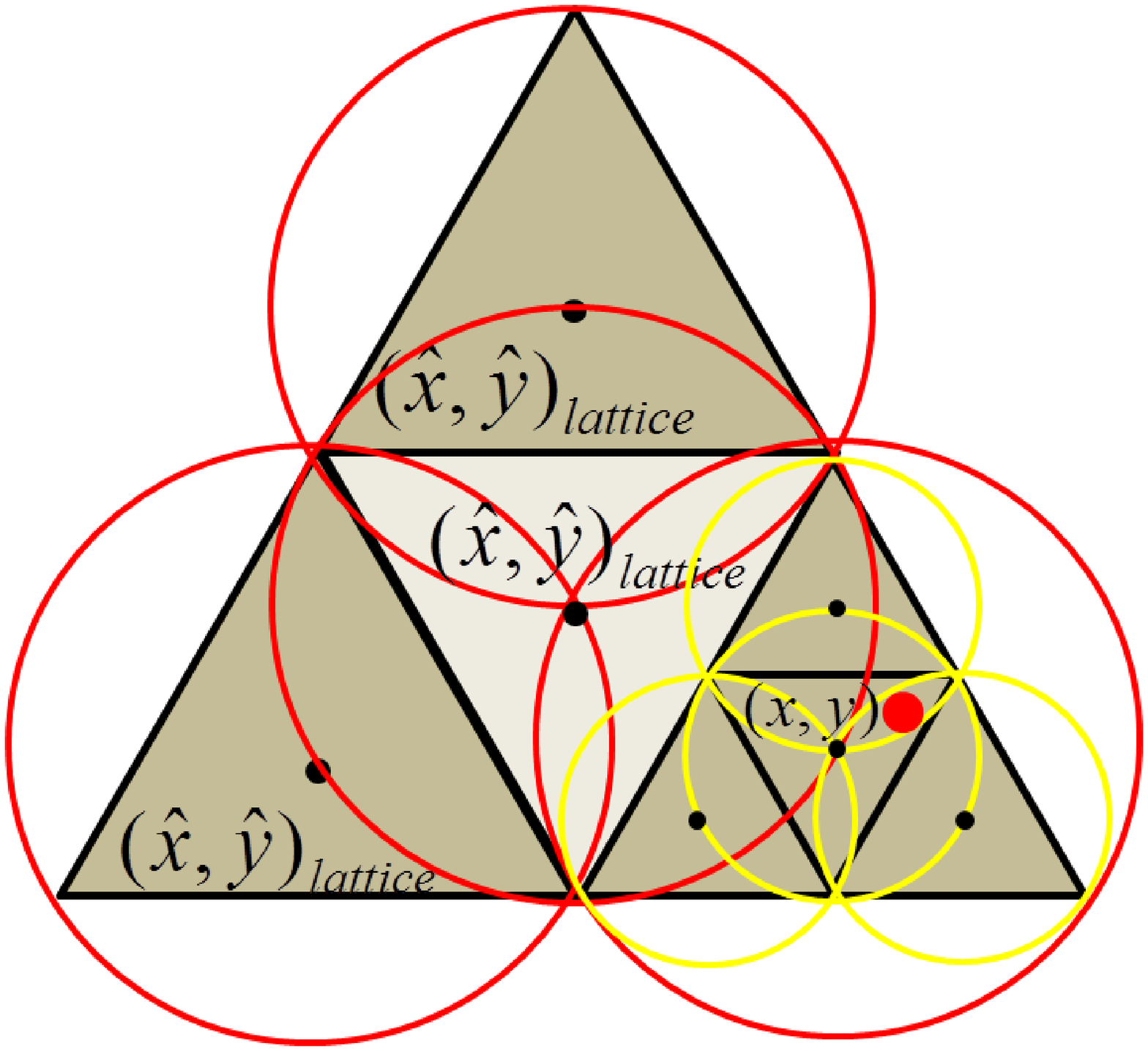}}
\caption{(a) Triangle orientation determined by cell index zero-parity. Translation magnitudes to adjacent neighbors scale by a factor of two per level. (b) Triangle edges are Voronoi diagrams of hexagonal lattice points at cell centers. A cell index is found via searching nearest lattice points within a current triangle.}
\end{center}
\end{figure}

To determine the orientation of index $t$ on level $l$, observe that the orientation only flips for center cells in Fig. \ref{fig:triagquadcenter}. This is the zero-parity of cell index $t$ (even is 0, odd is 1) adjusted for level $l$ and the number of components $|t|$. The function
\begin{equation}
\begin{split}
\textrm{isUp}(t, l) = 2 [ (\textrm{parity}(t) + l - |t| + 1) \bmod{2}] - 1,
\end{split}
\end{equation}
returns $1$ for the up orientation, and $-1$ for the inverted orientation. To determine the orientation of component $t_i$, make a query to $\textrm{isUp}(t_{1:i}, l-|t|+i)$ where cell index $t_{1:i}$ is the level adjusted ancestor. 
  
Note that triangle centers for components $t_i = 0$ are equivalent to its ancestor's center as no translations are necessary. The function $F(t, i)$ can now be expressed in polar coordinates
\begin{equation} 
\begin{split}
\displaystyle
j & =  l-|t|+i,  \quad 1 \leq t_i \leq 3, \\
\theta_i & = \textrm{isUp}(t_{1:i},j) \left ( \frac{\pi}{2} + \frac{2(t_i - 1) \pi }{3}  \right ), \quad R_i  =  2^{l_{max} - j} r,
\end{split}
\label{eq:triagcenter}
\end{equation}
where radius $R$ is scaled by a factor of two per level in Fig. \ref{fig:triagquadcenter}. The CellCenter($n$, $l$) procedure is as follows:
\begin{enumerate}
\item Let current center coordinates $(x,y) = (0,0)$,  $i=1$
\item Update center $(x,y) = (x,y) + (R_i \cos{\theta_i}, R_i \sin{\theta_i} )$ from eq. \ref{eq:triagcenter} \label{algo:triagcenter:2} 
\item Increment $i$ and repeat from step \ref{algo:triagcenter:2} until $i=l$
\item Return $(x,y)$
\end{enumerate}

\subsection*{CellIndex($x$, $y$, $l$)}

A cell index $t = n_7$ contains point $(x,y)$ on level $l$ if all ancestor cell indices $t_{1:l}$ also contain it. Since triangles descendants are contained within ancestor triangles, a search proceeds from level $0$ to $l$ s.t. successive descendant triangles contain point $(x,y)$.

Observe that triangle edges form Voronoi diagrams of regular hexagonal lattice points. That is, the triangle edges represent points equidistant between nearest lattice points. Hence, a query point $(x,y)$ that falls within a triangle is also the nearest neighbor of the triangle's center or lattice point in Fig. \ref{fig:triagquadindex}.

To track the progress of triangle centers, denote point $(u,v)$ as the current center. Note that triangle centers for components $t_i = 0$ are equivalent to its ancestor's center. That is, if a center triangle is found to be the nearest neighbor, then the current center remains unchanged while the vertical orientation $up$ flips. In polar coordinates, denote the three triangle lattice points about the origin as
\begin{equation} 
\begin{split}
\displaystyle
\theta_i & = up_i \left ( \frac{\pi}{2} + \frac{2(t_i - 1) \pi }{3}  \right ), \quad R_i  = 2^{l_{max} - i} r, \quad 1 \leq t_i \leq 3,
\end{split}
\label{eq:triaglattice}
\end{equation}
where similar to the CellCenter function, the radius $R$ is scaled by a factor of two per level. The CellIndex($x$, $y$, $l$) procedure is as follows:
\begin{enumerate}
\item Set to origin the current center $(u,v)= (0,0)$, let $i=1$
\item Assign component $t_i$ the cell type with the nearest lattice to query point $(x,y)$ on level $i$  \label{algo:triagindex:2} 
\item Update current center $(u,v) = (u,v) + (R_i \cos{\theta_i}, R_i \sin{\theta_i} )$ from eq. \ref{eq:triaglattice}
\item If $t_i = 0$, then flip orientation $up_i = -up_{i-1}$
\item Increment $i$ and repeat from step \ref{algo:triagindex:2} until $i=l$
\item Return $t$
\end{enumerate}

\section{Separation Ratios} \label{sec:separationratios}

\begin{figure}[ht]
\begin{center}
\subfloat[\label{fig:septreelocalsep}]{\includegraphics[scale=0.33]{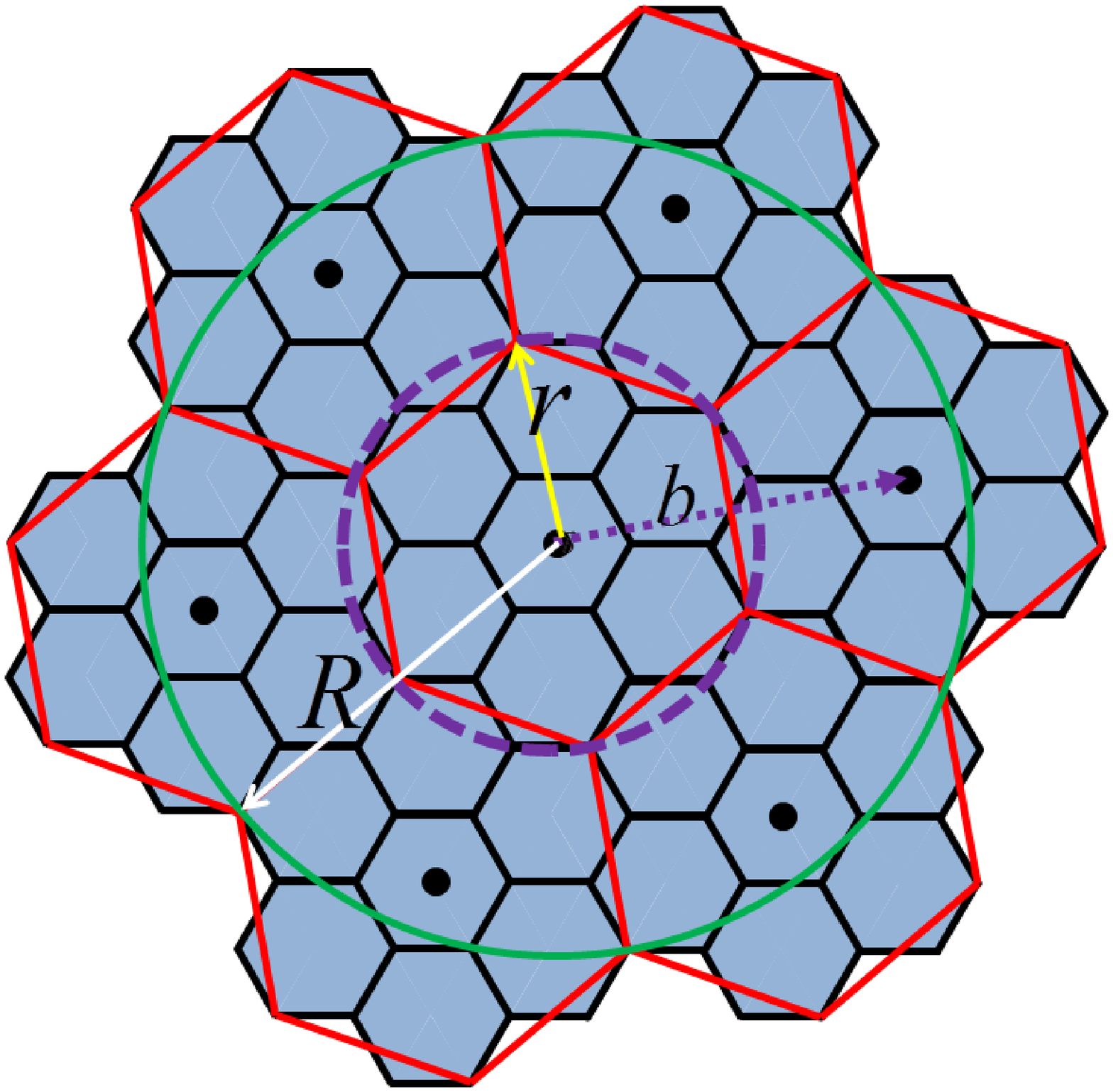}}\hspace{.3cm}
\subfloat[\label{fig:septreerho}]{\includegraphics[scale=0.35]{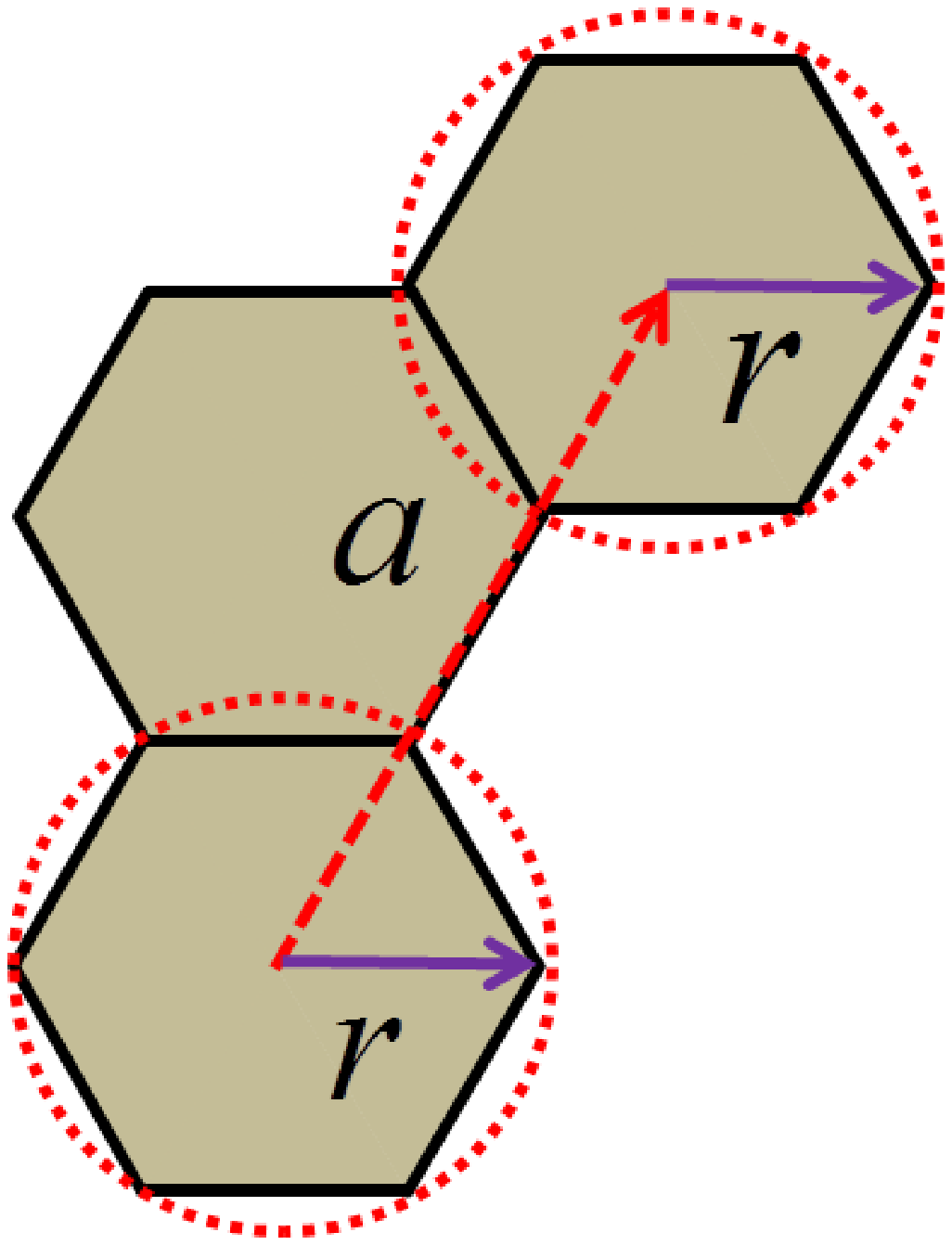}}
\caption{(a) Local separation ratio $r / R$ drawn from level $1$ lattice points and their induced hexagonal Voronoi diagrams. Either radius $r$ or $R$ is computable from magnitude $b$ obtained from eq. \ref{eq:triagcenter}. (b) WSPD distance $\rho = 3r-a = r$ generalized over level induced hexagonal boundaries}
\end{center}
\end{figure}

To show that septree local separation and WSPD ratios are level invariant, consider the base case in Fig. \ref{fig:hexagon} where ratio $r / R = 1/2$. For successive levels, lattice points induce Voronoi diagrams that form rotated and scaled hexagonal boundaries that intersect points along minor and major radii $r$ and $R$ in Fig. \ref{fig:septreelocalsep}. The rotation and scaling are proportional to eq. \ref{eq:triagcenter} adjusted for level $i$. The WSPD distance $\rho=r$ in Fig. \ref{fig:septreerho} and local separation ratio are preserved as the induced hexagonal boundaries are self-similar to the base case.

To show that the triangle-quadtree local separation ratio is level invariant, observe that ratio $r/R = 1/2$ for the same level in Fig. \ref{fig:triagquadcenter}. On successive levels, a major radius $R$ is equivalent to the minor radius $r$ from the preceding level and can be computed from eq. \ref{eq:triaglattice}. The WSPD ratio $\rho / r =\sqrt{7}-2$ in Fig. \ref{fig:rhotriag} is derived in appendix \ref{app:triagrho}. Similar to the quadtree, the level invariance stems from self-similar tiling arrangements with the base case in Fig. \ref{fig:triangle}.

\begin{figure}[ht]
\begin{center}
\subfloat[]{\includegraphics[scale=0.3]{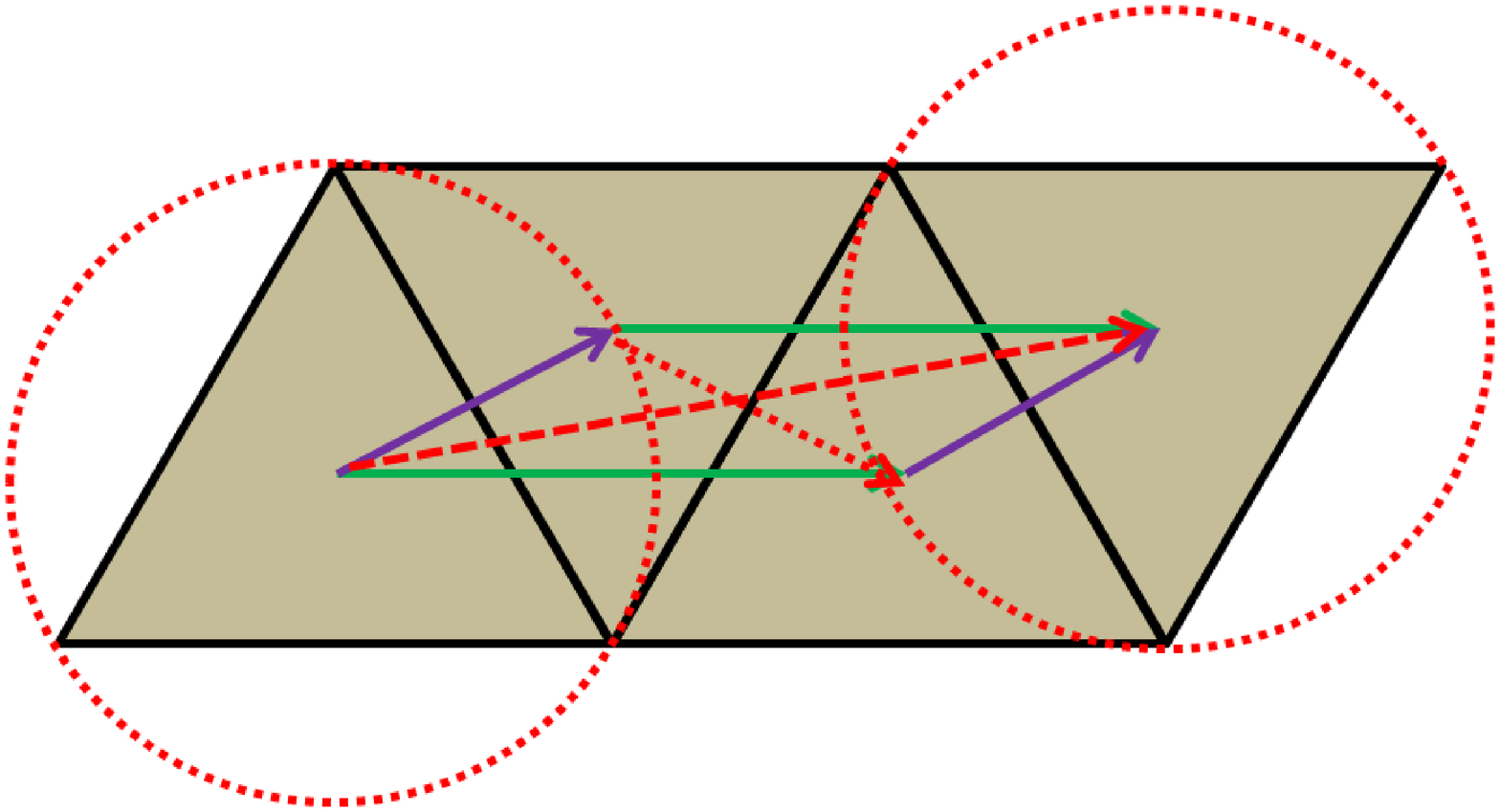}}\hspace{.3cm}
\subfloat[]{\includegraphics[scale=0.3]{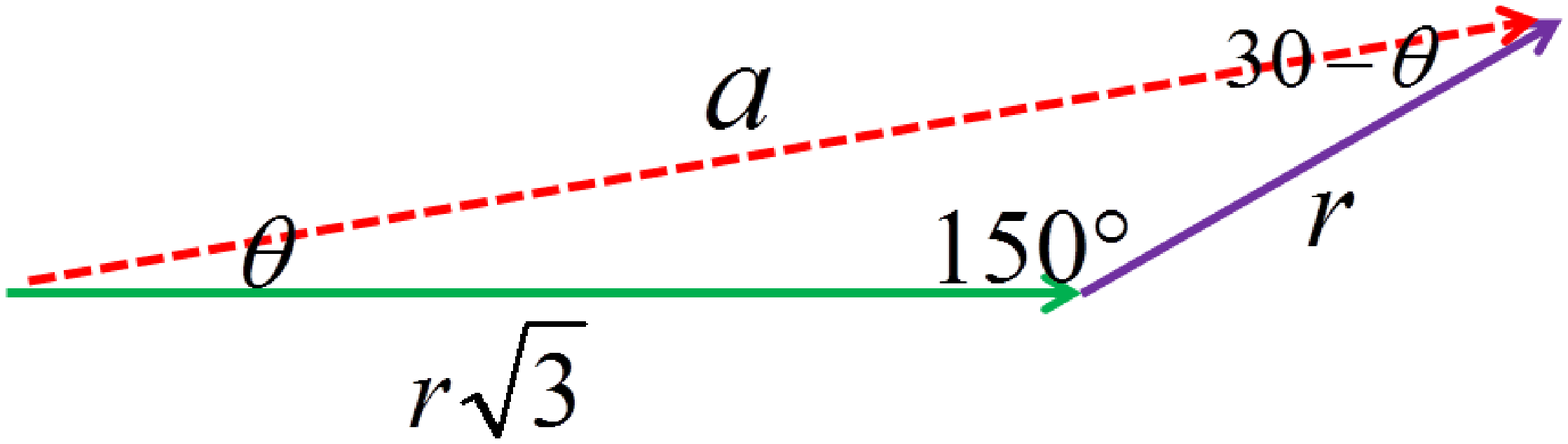}}
\caption{(a) Triangle tilings for nearest M2L translation. (b) Geometric analysis for WSPD distance } \label{fig:rhotriag}
\end{center}
\end{figure}

\section{Cost Analysis}\label{sec:prop}

To estimate the costs of a hierarchical FMM, a uniform source and target point distribution over the data structure's geometric domain is assumed \cite{FMM_COST}. The total computational costs w.r.t. the quadtree, septree and triangle quadtree are derived in appendices \ref{app:quad}, \ref{app:septree} and \ref{app:triagquad}. The overall costs depend on a density quantity $S_{opt}$ of source and target points per cell. This density quantity is a function of the number of M2L translations $P4$ and the number of neighbor cells $P2$ which is optimized by equating the number of M2L translations to the number of direct evaluations per cell. Substituting density $S_{opt}$ back into the FMM total costs yields the optimal FMM costs listed in table \ref{tab:props}.

The total FMM costs also depend on the number of truncation terms $p$ which are chosen a priori to guarantee a minimum error bound. This quantity depends on the kernel expansion and translation formulations as specified in \cite{COULOMBIC_MAST}. Maximum absolute error bounds \cite{COULOMBIC} are written in terms of local separation ratio $r/R$ and WSPD ratio $\rho / r$ between M2L translations. For a multipole expansion and translation, the absolute error is bounded by 
\begin{equation}
\begin{split}
\displaystyle
| \epsilon_p |  \leq \frac{\sum{|u_i|}_{i=1}^{k} }{R-r} \left ( \frac{r}{R} \right )^{p+1}.\\
\end{split}
\end{equation}
For a M2L translation, the absolute error is bounded by
\begin{equation}
\begin{split}
\displaystyle
| \epsilon_p | \leq \frac{\sum{|u_i|}_{i=1}^{k} }{\rho} \left ( \frac{1}{1+\rho / r} \right )^{p+1}.\\
\end{split}
\end{equation}
Last, the local to local translations are exact.

\begin{table}
\begin{center}
  \begin{tabular}{ | c | c | c | c | c | c | c |} 
    \hline
    \textbf{Struct} & \textbf{$P4$} & \textbf{$P2$} & \textbf{$r/R$} & \textbf{$\rho/r$} & \textbf{$S_{opt}$} & \textbf{Opt Costs} \\ \hline
    \textbf{Q-tree} & 27 & 9 & $\frac{\sqrt{2}}{3}$ & $2(\sqrt{2}-1)$ & $p \left (  \frac{116N}{27M} \right ) ^ {.5} $  & $ (M+N+32.64(MN)^{.5} )p$  \\ \hline
    \textbf{S-tree} & 42 & 7 & $\frac{1}{2}$ & $1$ & $p \left (  \frac{308N}{42M} \right ) ^ {.5} $ & $ (M+N+35.2(MN)^{.5} )p$ \\ \hline
    \textbf{TQ-tree} & 39 & 13 & $\frac{1}{2}$ & $\sqrt{7}-2$ & $p \left (  \frac{164N}{39M} \right ) ^ {.5} $ & $ (M+N+46.65(MN)^{.5} )p$ \\ \hline
  \end{tabular}    
  \caption{FMM properties for quadtree, septree and triangle quadtree where $p$ is the number of truncation terms for each expansion, P4 is number of M2L translations per cell, P2 number of cells in neighborhood, $\rho$ is the multipole to local separation distance and  $S_{opt}$ the optimal number of points per cell}
\label{tab:props}
\end{center}
\end{table}

For a fixed number of truncation terms, the total optimized cost for septree is slightly greater than that of quadtree as the lower neighbor count does not fully compensate for the greater number of M2L translations per cell. The triangle quadtree is less cost efficient in both categories and so yields a larger leading coefficient term. For multipole expansions, the larger local separation ratios $r/R$ in both septree and triangle quadtree suggest wider error bounds than that of quadtree. For M2L translations, the smaller WSPD ratio $\rho / r$ for septree suggests a tighter error bound than that of both quadtree and triangle-septree. This last property may have considerations in setting the maximum $l_{max}$ of the data structure.

\section{Experiments}\label{sec:experiments}

\begin{figure}[h]
\begin{center}
\includegraphics[scale=0.3]{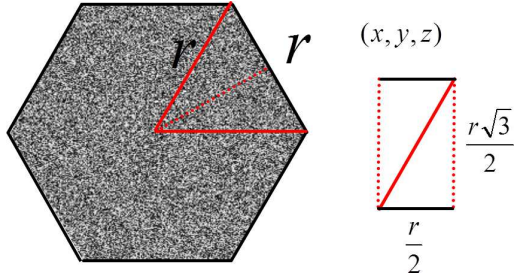}
\caption{Regular polygon point picking via parameterization into uniform variables $(x,y,z)$ where a point $(x,y)$ is uniformly sampled from a rectangle while uniformly random variable $z$ chooses the triangle}
\label{fig:uniform}
\end{center}
\end{figure}

To validate the theoretical FMM costs and error bounds from section \ref{sec:prop}, a set of test cases suitable for uniform conditions for each data structure is generated. That is, the input domain consist of $S_{opt}$ uniformly random source and target points fitted to each cell. A regular polygon point picking method in Fig. \ref{fig:uniform} handles the various tile geometries. For a regular cell, the polygon is subdivided into disjoint isosceles triangles per edge where halves of triangles are rearranged to form rectangles. Uniformly random points $(x,y)$ are picked from this rectangle and a third uniformly random variable $z$ assigns the point to one of the isosceles triangles. To modify the total number of source and target points, the total number of occupied cells is reduced.

For experiments, both runtime and error estimates of the FMM are considered while fixing either the number of source and target points ($N$, $M$), or the number of truncation terms $p$. Runtime comparisons w.r.t. the direct evaluation method are shown when appropriate. The implementation is written and tested on Matlab 2010b and Intel i7-2630QM hardware.

\begin{figure}[h]
\begin{center}
\includegraphics[scale=0.57]{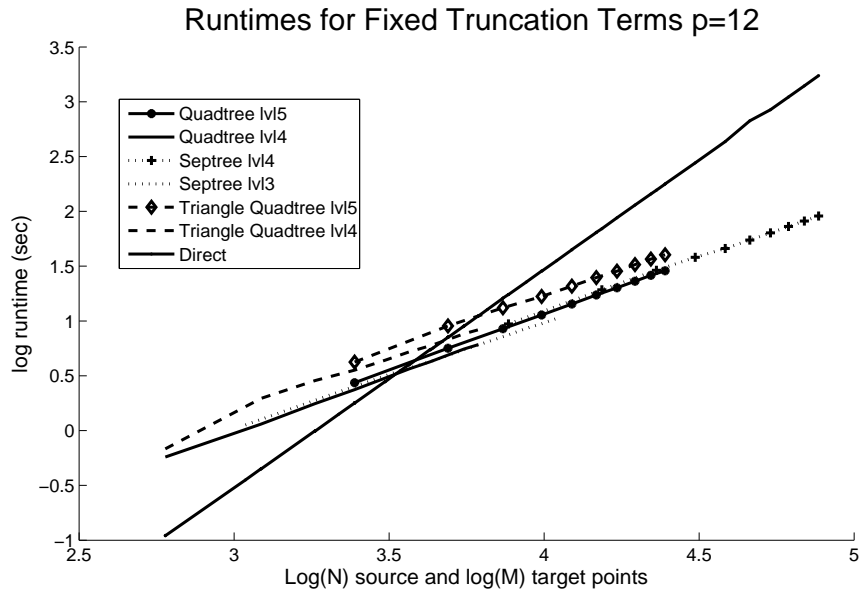}
\caption{Runtime (seconds) for FMM data structures with variable number of source $N$ and target $M$ points and fixed truncation terms $p=12$}
\label{fig:resrunn}
\end{center}
\end{figure}

\begin{figure}[h]
\begin{center}
\includegraphics[scale=0.57]{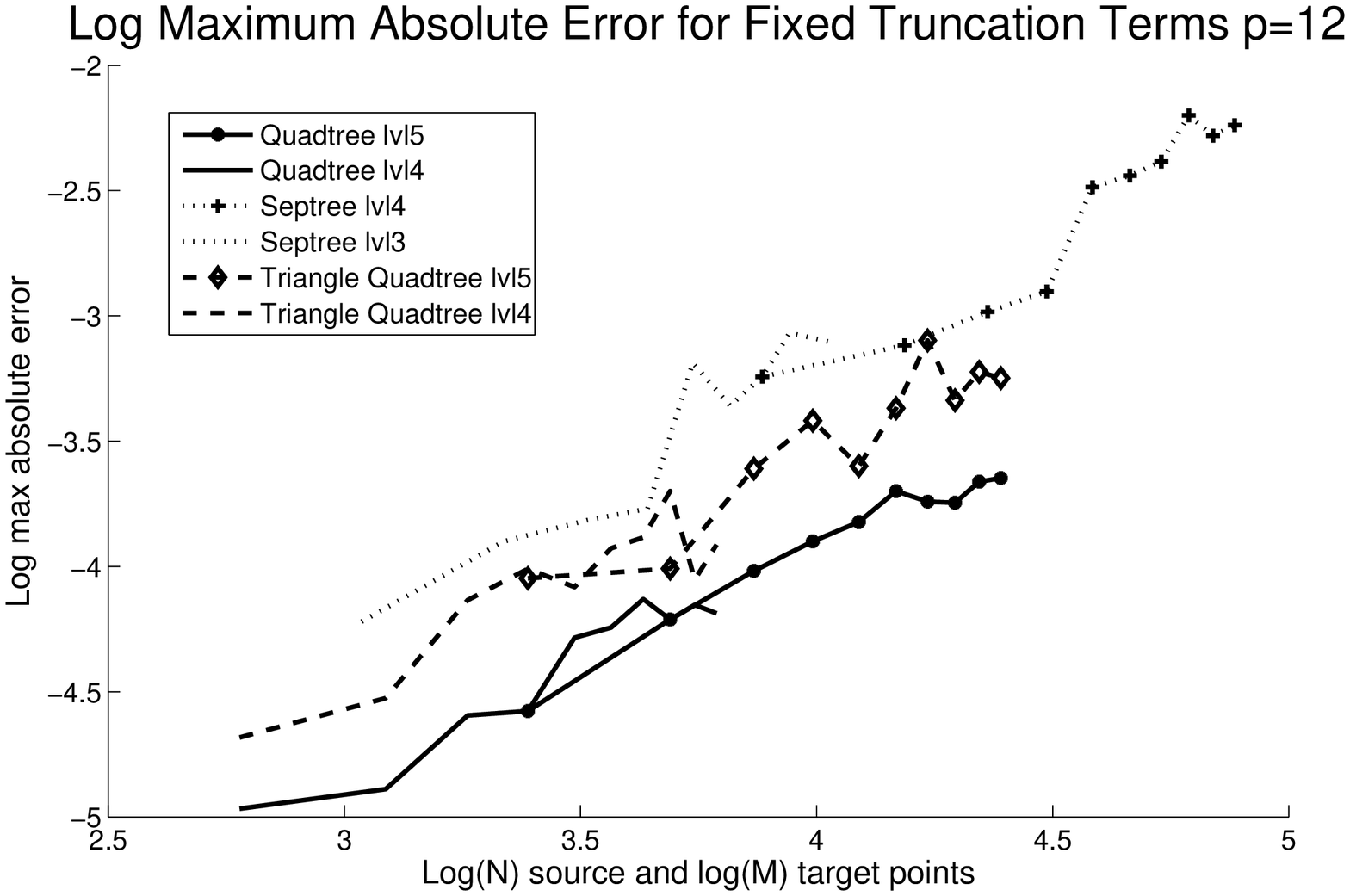}
\caption{Maximum absolute errors for FMM data structures with variable number of source $N$ and target $M$ points and fixed truncation terms $p=12$}
\label{fig:reserrn}
\end{center}
\end{figure}

The runtime results for a variable number of source and target points in Fig. \ref{fig:resrunn} match the theoretical estimates. All three hierarchical FMM data structures obtain linear asymptotic runtimes. The basic quadtree data structure outperforms the septree by a small margin for both max levels $l_{max}$. The triangle quadtree performs the slowest out of the three. The direct method obtains a quadratic runtime and where an estimated cross-over point between the number of source and target points is between $10^{3.5}$ to $10^{3.75}$. The error results in Fig. \ref{fig:reserrn} indicate a linear increase in error w.r.t. the number of source and target points $N$ and $M$. The quadtree obtains the least error presumably due to the smallest local separation ratio. One explanation for the greater maximum error in the septree compared to the triangle quadtree is the larger number of corner points in the tiling where source-target point distances are minimized. 

\begin{figure}[h]
\begin{center}
\includegraphics[scale=0.57]{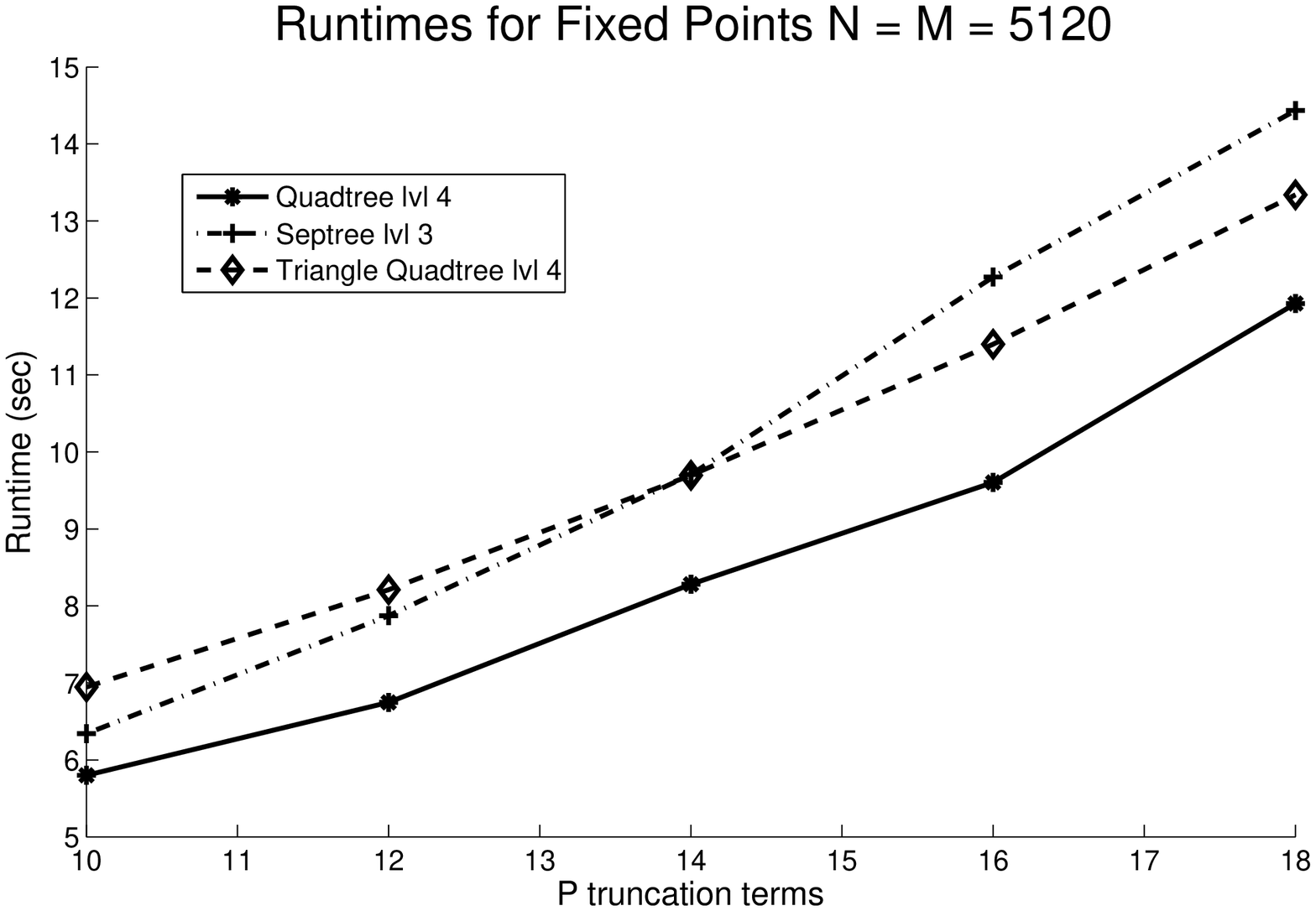}
\caption{Runtime (seconds) for FMM data structures with variable number of truncation terms $p$ and fixed number source and target points $N=M=5120$}
\label{fig:resrunp}
\end{center}
\end{figure}

\begin{figure}[h]
\begin{center}
\includegraphics[scale=0.57]{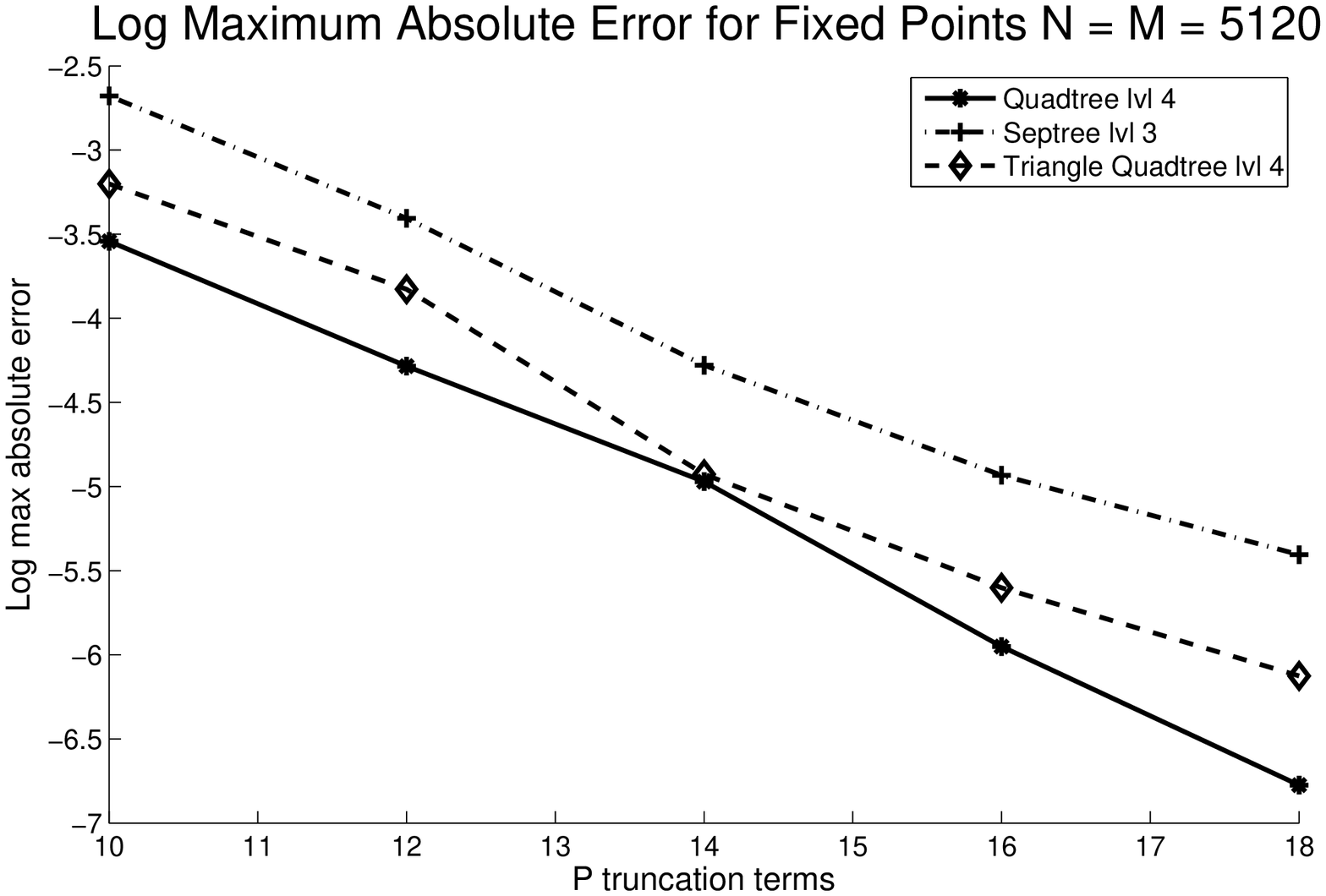}
\caption{Maximum absolute errors for FMM data structures with variable number of truncation terms $p$ and fixed number source and target points $N=M=5120$}
\label{fig:reserrp}
\end{center}
\end{figure}

The runtime results for a variable number of truncation terms $p$ in Fig. \ref{fig:resrunp} reveal a cross-over point between septree and triangle-quadtree. This may be due to the greater number of M2L translations in the septree and the inefficiencies in the implementation where the translation matrices are computed on the fly. The error results in Fig. \ref{fig:reserrp} exhibit an exponential error loss for increasing truncation terms.

\section{Conclusions} \label{sec:conclusions}
In this paper, we have shown that regular hexagonal and triangle tilings of the Euclidean plane generate septree and triangle-quadtree data structures that satisfy local separation and WSPD properties for the FMM. We derived their implementations and respective geometric properties with regards to FMM error bounds and costs. The empirical results validated the theoretical claims and have shown to be comparable to the original quadtree. While the quadtree remained ideal for a uniformly distributed data set, both the septree and triangle-quadtree may have applications for non-uniform domains, especially when the input is clustered about respective lattice points.

As a final remark, the septree's base hexagon unit is an order-3 permutohedron. Permutohedrons are $(d-1)$ dimensional polytopes embedded in $d$ dimensional space that can tessellate the domain and is indexable via a similar GBT addressing scheme. Future work may investigate constructions of a permutohedron-tree and its separation properties in higher dimension.

\section*{Acknowledgements}
We would like to thank Dr. Nail Gumerov for his lectures in the fast mutlipole methods course at the University of Maryland, College Park and partial support from the ONR Office of Naval Research under the MURI grant N00014-08-10638.

\bibliography{fmmbib}

\appendix

\section{FMM Quadtree Costs} \label{app:quad}
FMM evaluation costs before and after substituting optimal point density per cell for quadtree data structure in eqs. \ref{app:quadbefore} and \ref{app:quadafter}.
\begin{equation}
\displaystyle
\nonumber
\begin{split}
N = 2^{L_* d},\quad S = 2^{L_s d},\quad L = L_* - L_s, \quad K = 2^{dL}, \\ P_4(d) = 3^d(2^d-1), \quad P_2(d) = 3^d.
\end{split}
\end{equation}
%%%%%%%%%%%%%%%%%%%%%%%%%%%%%%%%%%%%%%%
\begin{equation}
\displaystyle
\begin{split}
\textrm{cost(FMM)} & = (M+N)P + \left ( K\frac{2^d}{2^d - 1} (P_4(d) + 2)-\frac{2^{3d+1}-P_4(d)2^{2d}}{2^d-1} \right ) P^2\\
		& +   M(P_2(d)s + P).
\end{split}\label{app:quadbefore}
\end{equation}
%%%%%%%%%%%%%%%%%%%%%%%%%%%%%%%%%%%%%%%
\begin{equation}
\displaystyle
\begin{split}
\textrm{opt cost} & = (M+N)P + \left ( \sqrt{\frac{2^d}{2^d - 1}} + \sqrt{\frac{2^d - 1}{2^d}} \right ) (MN (3^{d}(2^d-1)+2)3^d)^{.5}P \\
& \approx \left ( M+N +  32.64(MN)^{.5} \right ) P, \quad d=2.
\end{split}\label{app:quadafter}
\end{equation}

%%%%%%%%%%%%%%%%%%%%%%%%%%%%%%%%%%%%%%%
%%%%%%%%%%%%%%%%%%%%%%%%%%%%%%%%%%%%%%%
%%%%%%%%%%%%%%%%%%%%%%%%%%%%%%%%%%%%%%%
%%%%%%%%%%%%%%%%%%%%%%%%%%%%%%%%%%%%%%%
\section{FMM Septree Costs} \label{app:septree}
FMM evaluation costs before and after substituting optimal point density per cell for septree data structure in eqs. \ref{app:sepbefore} and \ref{app:sepafter}.
\begin{equation}
\displaystyle
\nonumber
\begin{split}
N = 7^{L_*}, \quad S = 7^{L_s}, \quad L =  L_* - L_s, \quad K = 7^{L}, \quad P_4 = 42, \quad P_2 = 7.
\end{split}
\end{equation}
%%%%%%%%%%%%%%%%%%%%%%%%%%%%%%%%%%%%%%%
\begin{equation}
\displaystyle
\begin{split}
\textrm{cost(FMM)} & = (M+N)P + \left ( K\frac{7}{6} (P_4 + 2)-\frac{7^{4}-P_47^{2}}{6} \right )P^2 +   M(P_2s + P).
\end{split}\label{app:sepbefore}
\end{equation}
%%%%%%%%%%%%%%%%%%%%%%%%%%%%%%%%%%%%%%%
\begin{equation}
\displaystyle
\begin{split}
\textrm{opt cost} & = (M+N)P + \left ( \sqrt{\frac{7}{6}} + \sqrt{\frac{6}{7}} \right ) (308MN)^{.5}P \\
		& \approx \left ( M+N + 35.2(MN)^{.5} \right ) P.
\end{split}\label{app:sepafter}
\end{equation}

%%%%%%%%%%%%%%%%%%%%%%%%%%%%%%%%%%%%%%%
%%%%%%%%%%%%%%%%%%%%%%%%%%%%%%%%%%%%%%%
%%%%%%%%%%%%%%%%%%%%%%%%%%%%%%%%%%%%%%%
%%%%%%%%%%%%%%%%%%%%%%%%%%%%%%%%%%%%%%%
\section{FMM Triangle Quadtree Costs} \label{app:triagquad}
FMM evaluation costs before and after substituting optimal point density per cell for triangle quadtree data structure in eqs. \ref{app:tquadbefore} and \ref{app:tquadafter}.
\begin{equation}
\displaystyle
\nonumber
\begin{split}
N = 4^{L_*},\quad S = 4^{L_s},\quad L = L_* - L_s, \quad K = 4^{L}, \quad P_4 = 39, \quad P_2 = 13.
\end{split}
\end{equation}
%%%%%%%%%%%%%%%%%%%%%%%%%%%%%%%%%%%%%%%
\begin{equation}
\displaystyle
\begin{split}
\textrm{cost(FMM)} & = (M+N)P + \left ( K\frac{4}{3} (P_4 + 2)-\frac{4^{4}-P_47^{2}}{3} \right )P^2 +   M(P_2s + P).
\end{split} \label{app:tquadbefore}
\end{equation}
%%%%%%%%%%%%%%%%%%%%%%%%%%%%%%%%%%%%%%%
\begin{equation}
\displaystyle
\begin{split}
\textrm{opt cost} & = (M+N)P + \left ( \sqrt{\frac{4}{3}} + \sqrt{\frac{3}{4}} \right ) (533MN)^{.5}P \\
		& \approx \left ( M+N + 46.65(MN)^{.5} \right ) P.
\end{split} \label{app:tquadafter}
\end{equation}

%%%%%%%%%%%%%%%%%%%%%%%%%%%%%%%%%%%%%%%
%%%%%%%%%%%%%%%%%%%%%%%%%%%%%%%%%%%%%%%
%%%%%%%%%%%%%%%%%%%%%%%%%%%%%%%%%%%%%%%
%%%%%%%%%%%%%%%%%%%%%%%%%%%%%%%%%%%%%%%
\section{GBT Translations} \label{app:Septreetrans}
Septree GBT translations for zero-extending indices.
The angular separation $\theta$ in Fig. \ref{fig:GBTextendedgeometry} is
\begin{equation}
\displaystyle
\begin{split}
\frac{\sin{\theta}}{2a} = \frac{\sin{120}}{b} = \frac{\sin{(60-\theta)}}{a} \Rightarrow \theta = \arctan{\frac{\sqrt{3}}{2}}.
\end{split}\label{APP:GBTANG}
\end{equation}

The magnitude of translation vector $\vec{b}$ in Fig. \ref{fig:GBTextendedgeometry} is
\begin{equation}
\displaystyle
\begin{split}
b\sin{\arctan{\frac{\sqrt{3}}{2}} } = b \frac{\frac{\sqrt{3}}{2}}{\sqrt{1+ \left ( \frac{\sqrt{3}}{2} \right )^2 }} = a\sqrt{3} \Rightarrow  b = a \sqrt{7}, \quad a_0 = r\sqrt{3}.
\end{split}\label{APP:GBTMAG}
\end{equation}

%%%%%%%%%%%%%%%%%%%%%%%%%%%%%%%%%%%%%%%
%%%%%%%%%%%%%%%%%%%%%%%%%%%%%%%%%%%%%%%
%%%%%%%%%%%%%%%%%%%%%%%%%%%%%%%%%%%%%%%
%%%%%%%%%%%%%%%%%%%%%%%%%%%%%%%%%%%%%%%
\section{Triangle-quadtree WSPD} \label{app:triagrho}
The min M2L translation distance $\rho$ for triangle quadtree in Fig. \ref{fig:rhotriag} is
\begin{equation}
\displaystyle
\begin{split}
\frac{\sin{\theta}}{r} = \frac{\sin{150}}{a} = \frac{30-\theta}{r\sqrt{3}} \Rightarrow \theta = \arctan{\frac{\sqrt{3}}{9} }, \\
\frac{r}{2} = a\sin{\arctan{\frac{\sqrt{3}}{9}} } = a \frac{\frac{\sqrt{3}}{9}}{\sqrt{1+ \left ( \frac{\sqrt{3}}{9} \right )^2 }} = a \frac{1}{2\sqrt{7}}, \\
a = r \sqrt{7}, \quad \rho = r(\sqrt{7}-2).
\end{split}
\end{equation}

\end{document}